\documentclass{amsart}
\usepackage{amssymb,amscd}
\usepackage[arrow,matrix]{xy}
\newcommand\datver[1]{\def\datverp% 
{\par\boxed{\boxed{\text{Version: #1; Run: \today}}}}} 
\datver{3.0A; Revised: August 14, 2001}

\newcommand\Hc{\operatorname{HC}}
\newcommand\Hp{\operatorname{HP}}
\newcommand\tHp{\operatorname{HP}^{\topo}}
\newcommand\Hd{\operatorname{HH}}
\newcommand\tHd{\operatorname{HH}^{\topo}}
\newcommand\topo{\operatorname{top}}
\newcommand\opp{\operatorname{op}}
\newcommand\mfk{\mathfrak}

\newcommand\pa{\partial}

\newcommand\supp{\operatorname{supp}}
\newcommand\cohom{\operatorname{H}}

\newcommand\CC{\mathbb C}

\newcommand\FF{\mathbb F}
\newcommand\KK{\mathbb K}

\newcommand\kk{\mathbf k}
\newcommand\ckk{\widehat{\kk}}

\newcommand\QQ{\mathbb Q}
\newcommand\ZZ{\mathbb Z}
\newcommand\Wh{\widehat}
\newtheorem{theorem}{Theorem}
\newtheorem{proposition}{Proposition}
\newtheorem{corollary}{Corollary}
\newtheorem{lemma}{Lemma}
\newtheorem{definition}{Definition}
\theoremstyle{remark}

\newcommand\ie{{\em i.e., }}
\newcommand\Prim{\operatorname{Prim}}
\newcommand\prn{\operatorname{Prim}_n(A)}
\newcommand\PR[1]{\operatorname{Prim}(#1)}
\newcommand\prim[2]{\operatorname{Prim}_{#1}(#2)}
\newcommand\Jac{\operatorname{Jac}}

\newcommand\End[2]{\operatorname{End}_{#1}(#2)}

\newcommand\Tt{a_0 \otimes a_1 \otimes \ldots \otimes a_n}

\newcommand\Hecke{{H}}
\newcommand\gHecke{{\mathbb{H}}}
\newcommand\param{{\bf \overline{q}}}
\newcommand\pp{\mfk p}
\newcommand\Max{\operatorname{Max}}
\newcommand\hW{\widehat W}

\begin{document}

\title[Homology of Hecke algebras]
{Periodic cyclic homology of Iwahori-Hecke algebras}

\author[P. Baum]{Paul Baum} \address{Pennsylvania State University,
	 Math. Dept., University Park, PA 16802, USA}
	 \email{baum@math.psu.edu} \author[V. Nistor]{Victor Nistor}
	 \address{Pennsylvania State University, Math. Dept.,
	 University Park, PA 16802, USA} \email{nistor@math.psu.edu}
	 \thanks{P.B. was partially supported by NSF Grant
	 DMS-9704001.  V.N. was partially supported by the NSF Young
	 Investigator Award DMS-9457859, NSF Grant DMS-9971951 and
	 ``collaborative CNRS/NSF research grant'' DMS-9981251.
	 Manuscript available from {\bf \small www.math.psu.edu/baum/}
	 or {\bf \small www.math.psu.edu/nistor/}}
	 
\dedicatory\datverp

\begin{abstract}\ 
We determine the periodic cyclic homology of the Iwahori-Hecke
algebras $\Hecke_q$, for $q \in \CC^*$ not a ``proper root of unity.''
(In this paper, by a {\em proper root of unity} we shall mean a root
of unity other than 1.) Our method is based on a general result on
periodic cyclic homology, which states that a ``weakly spectrum
preserving'' morphism of finite type algebras induces an isomorphism
in periodic cyclic homology. The concept of a weakly spectrum
preserving morphism is defined in this paper, and most of our work is
devoted to understanding this class of morphisms. Results of
Kazhdan--Lusztig and Lusztig show that, for the indicated values of
$q$, there exists a weakly spectrum preserving morphism $\phi_q :
\Hecke_q \to J$, to a fixed finite type algebra $J$. This proves that
$\phi_q$ induces an isomorphism in periodic cyclic homology and, in
particular, that all algebras $\Hecke_q$ have the same periodic cyclic
homology, for the indicated values of $q$. The periodic cyclic
homology groups of the algebra $\Hecke_1$ can then be determined
directly, using results of Karoubi and Burghelea, because it is the
group algebra of an extended affine Weyl group.
\end{abstract}

\maketitle
 
\tableofcontents
\section*{Introduction}

Let $G$ be the set of rational points of a reductive group defined
over a $p$-adic field $\FF$.  Fix an Iwahori subgroup $I \subset
G$. The space of complex, compactly supported functions on $G$ that
are right and left invariant with respect to $I$, with the convolution
product, is called the {\em Iwahori-Hecke algebra of $G$}. Since the
Iwahori subgroups of $G$ are all conjugate, the isomorphism class of
the Iwahori-Hecke algebra of $G$ is independent of the choice of the
Iwahori subgroup $I$.

The Iwahori-Hecke algebras $\Hecke_q$ were introduced in
\cite{IwahoriMatsumoto} by Iwahori and Matsumoto, who also provided in
that paper a description of these algebras using generators and
relations. Let $\hW$ be the extended affine Weyl group of $G$. As a
vector space $\Hecke_q \cong \CC[\hW]$, which provides us with a
distinguished system of generators of $\Hecke_q$. We shall denote the
elements of this system of generators by $T_w$, $w \in \hW$.  The
relations among the generators $T_w$ of the Iwahori-Hecke algebra
$\Hecke_q$, depend on the parameter $q$. By allowing $q$ to be an
arbitrary non-zero complex number, one obtains a family of
algebras $\Hecke_q$, which specialize to the Iwahori-Hecke algebra of
$G$ when $q$ is the number of elements in the residue field of $\FF$.
(See Section~\ref{Sec.Hecke}, where we review these definitions and
constructions in detail.)

Iwahori-Hecke algebras play a central role in the representation
theory of $p$-adic group. For example, the classification of the
irreducible, smooth representations of $G$ with non-zero vectors fixed
by $I$ reduces to the classification of the irreducible
representations of the Iwahori-Hecke algebra associated to $G$ (see
\cite{Borel}). Other representations of $G$ can also be classified in
terms of Iwahori-Hecke algebras. Indeed, if $G = GL_n(\FF)$, the
general linear group over the $p$-adic field $\FF$, all irreducible,
smooth representations of $G$ are classified by representations of
(finite tensor products of) Iwahori-Hecke algebras~\cite{BushnellKutzko}.

The representations of Iwahori-Hecke algebras (for suitable values of
$q$ and suitable groups $G$) were classified by Kazhdan and Lusztig in
a ground breaking paper, \cite{KL}, using the $K$-homology of various
varieties of Borel subgroups of the Langlands dual of $G$.  The group
$G$ does not appear in the Kazhdan-Lusztig classification, which is
based on the algebraic description of the Iwahori-Hecke algebras in
terms of generators and relations. However, the study of the
representation theory of $G$ was certainly their main motivation.

In this paper, by a {\em proper root of unity} we shall mean a root of
unity other than~1. Our work gives a complete determination of the
periodic cyclic homology of the Iwahori-Hecke algebras $\Hecke_q$,
provided that $q \in \CC^*$ is not a proper root of unity.  The first
step of our computation is to show that the periodic cyclic homology
groups $\Hp_*(\Hecke_q)$ are naturally isomorphic for all these values
of $q$. For $q = 1$ the Iwahori-Hecke algebra, $H_1$, is
the group algebra of an extended affine Weyl group, so we can use the
results of Karoubi \cite{Karoubi} and Burghelea \cite{Burghelea} to
determine the groups $\Hp_*(\Hecke_1)$ directly. This leads to a
complete and natural determination of the groups $\Hp_*(\Hecke_q)$,
for $q \in \CC^*$ not a root of unity. (The approach in \cite{BC} is
in fact more convenient for the calculation, and we shall discuss
this in \cite{BBN}.) An important concept that facilitates our
calculations is that of finite type algebra, 
Definition~\ref{Def.finite.type}. Actually, a crucial step in our
calculations is the extension of a technical result on finite type
algebras from \cite{KNS} (see Theorem~\ref{theorem.Mult}).

The idea of the proof is to study a certain class of morphisms of
finite type algebras that we introduce in this paper, namely the class
of ``weakly spectrum preserving morphisms,''
Definition~\ref{def.w.sp.p}. The relevance of this class is
two-fold. First, a spectrum preserving morphism defines an isomorphism
in periodic cyclic homology (this is the content of
Theorem~\ref{theorem.HPiso}), which is one of the main results of this
paper. Second, Lusztig has constructed an algebra $J$, independent of
$q$, and morphisms
\begin{equation}
	\phi_q : \Hecke_q \to J, 
\end{equation}
which he used to classify the irreducible representations of
$\Hecke_q$ in terms of the irreducible representations of $J$.  One of
the main steps in this classification is a result, which, in our
language, amounts to saying that $\phi_q$ is weakly spectrum
preserving for all values of $q \in \CC^*$ that are not proper roots
of unity.  See \cite{KL,LC1,LC2,LC3,LC4}, and especially
\cite{Asterisque}.  Thus, by Theorem~\ref{theorem.HPiso},
\begin{equation}
	(\phi_q)_* : \Hp_*(\Hecke_q) \to \Hp_*(J)
\end{equation}
is an isomorphism, provided that $q \in \CC^*$ is not a proper root of
unity. Consequently, all groups $\Hp_*(\Hecke_q)$ are isomorphic for
our values of $q$, as mentioned also above. Moreover, these
isomorphisms are canonical.

The group $G$ will not appear in the statements or the proofs of our
results. Nevertheless, as in the Kazhdan-Lusztig theory, our main
motivation is to study the representation theory of $G$. Let us
briefly explain how our results could be related to the representation
theory of $G$.  Let $\mathcal{H}(G)$ be the (full) Hecke algebra of
$G$, \ie the algebra of compactly supported, locally constant
functions on $G$ with the convolution product being defined using a
fixed Haar measure. Then the groups $\Hp_*(\mathcal{H}(G))$ have two
different descriptions: one using the geometric (and algebraic)
structure of $G$ and another one using the representation theory of
$G$.

The first, algebraic, description of the groups $\Hp_*({\mathcal
H}(G))$ was obtained independently in \cite{Higson-Nistor1,Schneider1}
using the action of $G$ on its affine Bruhat-Tits building, and thus not 
relying on the representation theory of $G$. An even more precise, geometric
description of $\Hp_*({\mathcal H}(G))$ in terms of certain data that
carry arithmetic information was obtained in \cite{NistorHOMPDC}. This
description is in terms of conjugacy classes of elements of $G$ and
the continuous cohomology of the stabilizers of the semi-simple
conjugacy classes. See also \cite{BHP-hc}.

Let $\Wh{G}$ be the admissible dual of $G$, that is, the space of
smooth, irreducible representations of $G$ with the Jacobson topology.
In \cite{KNS}, the groups $\Hp_*({\mathcal H}(G))$ were related to the
de Rham cohomology of the strata of $\Wh{G}$ via a spectral
sequence.  This provides a representation theoretic, or spectral,
description of $\Hp_*({\mathcal H}(G))$. Relating the geometric and
spectral descriptions of the periodic cyclic cohomology groups of
$\mathcal H(G)$ will thus translate into concrete relations between
the geometric and spectral properties of $G$.

Let $D$ be a connected (or Bernstein) component of
$\Wh{G}$.  As a topological space, each Bernstein component $D$
is homeomorphic to the primitive ideal spectrum of a finite type
algebra $A_D$, canonically associated to $D$, and hence we have the
following disjoint union decomposition of $\Wh G$:
\begin{equation} \label{eq.red1}
	\Wh G = \cup_D\PR{A_D}.
\end{equation}
The corresponding relation in periodic cyclic homology is that
\begin{equation} \label{eq.red2}
	\Hp_*({\mathcal{H}(G)}) \cong \oplus_D \Hp_*(A_D).
\end{equation}
The proof of \eqref{eq.red2} uses the fact that the Hochschild
homology of the full Hecke algebra of $G$ vanishes above the split
rank of $G$. (There exist finite type algebras whose Hochschild
homology is non-zero in infinitely many dimensions.)

One can formulate then a strategy of relating
$\Hp_*({\mathcal{H}(G)})$ to the topology of $\Wh G$ as
follows. First, using \eqref{eq.red2}, we reduce the computation of
$\Hp_*({\mathcal{H}(G)})$ to the computation of $\Hp_*(A_D)$. Next,
the periodic cyclic homology groups of the finite type algebras $A_D$ are
related to the topology of their primitive ideal spectrum by the results
of \cite{KNS} and of this paper.  Then, finally, we determine the
topology of $\Wh G$ from that of the individual Bernstein
components using \eqref{eq.red1}. For this strategy to be effective,
we need a determination of the groups $\Hp_*(A_D)$ that is as explicit
as possible. Our results provide complete results for the Bernstein
components whose commuting algebras are Iwahori-Hecke algebras (or
tensor products of such algebras).

Another motivation to study the periodic cyclic homology of
Iwahori-Hecke algebras is that, as explained by Baum, Higson, and
Plymen in \cite{BHP}, the calculation of the periodic cyclic homology of Iwahori-Hecke
algebras can be used to prove the Baum-Connes conjecture for $GL_n(\FF)$. (See
\cite{BCH} for a statement of the Baum-Connes conjecture for $p$-adic
groups). See also \cite{BHP2}. It would be interesting to relate our
approach to the computation of the periodic cyclic homology of Iwahori-Hecke
algebras with the description of Iwahori-Hecke algebras in
\cite{GinzburgKV}. The implicit deformation idea used in this paper
fits also with the philosophy of quantum groups \cite{Manin,LusztigQG}.

We thank Alain Connes, Edward Formanek, Nigel Higson, David Kazhdan,
George Lusztig, Roger Plymen, and Peter Schneider for useful
discussions. The second named author would also like to thank 
the Max Planck Institut f\"ur Mathematik in Bonn for support
and hospitality while parts of this work have been completed.

%%%%%%%%%%%%%%%%%%%%%%%%%%%%%%%%%%%%%%%%%%%%%%%
%%%%%%%%%%%%%%%%%%%%%%%%%%%%%%%%%%%%%%%%%%%%%%%

\section{Finite type algebras\label{Sec1}}

In this section, we recall some facts about finite type algebras and 
prove a few needed properties of finite type algebras. We also fix notation.

In what follows, $Z(\mfk A)$ will denote the center of an  algebra
$\mfk A$. An algebra with unit will be called a {\em unital}
algebra. By $\kk$ we shall denote a commutative noetherian ring (with
unit), which in most cases will be a finitely generated commutative
algebra.  By a $\kk$--algebra (without unit in general), we shall mean
an algebra $A$ which is a $\kk$--module such that its algebraic
operations are $\kk$-linear. If $A$ is unital, this is the same as
saying that there is given a unital morphism $\kk \to Z(A)$.

For any $\kk$-algebra $A$, we shall denote by $A^+:=A \oplus \kk$, the
algebra obtained from $A$ by adjoining a (possibly new) unit.  The
multiplication is given by the straight-forward formula:
\begin{equation*}
	(a,\lambda)(b,\mu) = (ab + \lambda b + \mu a, \lambda \mu).
\end{equation*}
A $\kk$-linear morphism $A \to B$ of two $\kk$-algebras extends in an
obvious way to a unital $\kk$-linear morphism $A^+ \to B^+$ of the
algebras with adjoined unit.

\begin{definition}\label{Def.finite.type} \
Let $\kk\cong \CC[X_1,\ldots,X_r]/I$ be a finitely generated complex
ring.  A {\em finite type $\kk$--algebra} is a (not necessarily
unital) $\kk$--algebra $A$ that is a finitely generated
$\kk$--module. A {\em finite type algebra} is an algebra that is a
finite type $\kk$--algebra for some $\kk$.
\end{definition}

Thus while in this paper $\kk$ may occasionally denote a more general
commutative ring, when we are talking of a finite type $\kk$-algebra,
we always require $\kk$ to be a finitely generated complex ring. 

Recall that a primitive ideal of $A$ is a two-sided ideal $\mfk P
\subset A$ that is the kernel of a {\em non-zero} irreducible complex
representation of $A$ (so $\mfk P \not = A$).  Denote by $\PR A$ the
primitive ideal spectrum of $A$ endowed with the Jacobson
topology. Recall that the closed sets of the Jacobson topology are the
sets of the form $V(S_0)$,
\begin{equation*}
	V(S_0) := \{ \mfk P \in \PR{A},\ S_0 \subset \mfk P \,\}  ,
\end{equation*}
$S_0 \subset A$ arbitrary. Clearly, $V(I) = \Prim(A/I)$, for any
two-sided ideal $I$ of $A$. Also, recall that the {\em Jacobson
radical} of $A$, $\Jac(A)$, is the intersection of all primitive
ideals of $A$. If $\Jac(A) = 0$, we say that $A$ is {\em
semi-primitive}.

We shall use several times the following well known result (see
\cite{Rowen1} or \cite{KNS}, Lemma~1).

\begin{lemma}\label{lemma.onto} \ 
Let $A$ be a unital finite type algebra with center $Z=Z(A)$, and let
$\mfk P$ be a primitive ideal of $A$, then $\Theta(\mfk P) := \mfk P
\cap Z$ is a maximal ideal of $Z$ and $A/\mfk P \cong M_k(\CC)$ for
some $k$. If $A$ is generated by $q$ elements as a $Z$-module, then
$k^2 \le q$.
\end{lemma}

\begin{proof} \ 
Since $Z$ has an at most countable dimension as a complex vector
space, it follows that $A$ and any simple left $A$--module $V$ also
have an at most countable $\CC$-dimension.  The center $Z_0$ of $\End
A V$ is a field extension of $\CC$, because $V$ is irreducible. Since
the only field extension of finite or countable dimension of $\CC$ is
$\CC$ itself, we obtain that $Z \to Z_0$ is surjective, and hence
$\mfk P \cap Z$ is a maximal ideal of $Z$.  Moreover, $ Z/ Z \cap
\mfk P \cong \CC$, (Hilbert's Nullstellensatz).

We obtain then that, $A/\mfk P \cong M_n(\CC)$, because all centrally simple
complex algebras are matrix algebras.  If $A$ is generated as a
$Z$--module by $q$ elements, then $A/ \mfk P$ is generated by $q$
elements as a $Z/ Z \cap \mfk P$-module, so $n \le q^2$.
\end{proof}

See also \cite{LC3} and \cite{Quillen}.

Let $A$ be a finite type $\kk$-algebra.  If $\mfk P \subset A$ is the
kernel of an irreducible representation $(\pi,V_\pi)$, then the
isomorphism class of $\pi$ and the dimension of the space $V_\pi$ on
which it acts are uniquely determined. This explains why studying the
set $\Prim(A)$ of primitive ideals of $A$ is the same thing as
studying the irreducible representations of $A$.

We shall denote by $d_{\mfk P}:=\dim V_\pi$ the dimension of the space
$V_\pi$ on which $\pi$ acts and by $\prim{n}{A}$ the subset of
$\Prim(A)$ consisting of those primitive ideals $\mfk P$ for which
$d_{\mfk P} = n$.  We know that $d_{\mfk P}$ is not greater than the
number of generators of $A$ as a $\kk$--module, for all $\mfk P \in
\Prim(A)$, and hence $\prim{n}{A}$ is empty if $n$ is greater than a
fixed number $N$ that depends on $A$.

We shall denote by $\Max(\kk)$ the maximal ideal spectrum of a
commutative ring $\kk$ ($\Max(\kk) = \Prim(\kk)$ if $\kk$ is finitely
generated). Let $A$ be a finite type $\kk$-algebra. Then
Lemma \ref{lemma.onto} gives rise to a map
\begin{equation}
\label{inf.character}
	\Theta : \PR A \ni {\mathfrak P} \longrightarrow {\mathfrak p}=
	Z(A) \cap {\mathfrak P}\in \Max(Z(A)) = \Prim(Z(A)),
\end{equation}
called the {\em central character}. We can then define in a similar
way a map
\begin{equation}\label{inf.character2}
	\Theta_A : \PR A \ni {\mathfrak P} \longrightarrow {\mathfrak
	p} := \sqrt{\kk \cap {\mathfrak P}} = \kk \cap {\mathfrak P}
	\in \Max(\kk),
\end{equation}
also referred to as central character. The central character $\Theta_A$
map satisfies
\begin{equation*}
	\Theta_A^{-1}(V(I)) = V(IA),
\end{equation*}
for any ideal $I \subset \kk$, and hence it is continuous and surjective (see
\cite{KNS}, Lemma~1).

Every irreducible representation $\pi$ of a finite type $\kk$-algebra
$A$ extends uniquely to an irreducible representation $\pi^+$ of
$A^+$. If $\mfk P$ is the kernel of $\pi$, then $a \in \kk$ acts on
the space of $\pi$ via the scalar $\pi^+(a)$.  In particular, this
provides us with a way of defining the central character for
non-unital algebras.

We now discuss the important and well known result on the nilpotency
of the Jacobson radical of finite type algebras, Lemma
\ref{lemma.nilp}. We shall use this result several times in what
follows, so we include a proof, which is simpler in our setting of
finite type algebras then in the most general case. See
\cite{Braun,Kemer} and \cite{Procesi2} for the general case. The ideas
of our proof will be useful again later on. A different proof of this
result for finite type algebras was also given in \cite{KNS} (note
however that in that proof one has to first replace $A$ with
$\Jac(A)^+$, which is always possible).

If $\pp \subset \kk$ is a prime ideal and $M$ is a $\kk$-module, then
we shall denote by $M_\pp$, as is customary, the {\em localization
of $M$ at $\pp$}, that is,
\begin{equation*}
	M_\pp = S^{-1}M,
\end{equation*}
where $S = \kk \smallsetminus \pp$. Recall then that the support of a
$\kk$-module $M$ (a subset of $\Max(\kk)$) is defined as the set of
maximal ideals $\pp \subset \kk$ such that the localization $M_\pp
\not = 0$. Also, recall that the {\em dimension} of $\kk$ is the
maximum length of chains of prime ideals in $\kk$, see
\cite{AtiyahMacDonald1}, Chapter 8. See also \cite{Harris}. A chain
consisting of $d + 1$ distinct ideals has length $d$. Non-proper
ideals, that is $(0)$ and $\kk$ are allowed in such a chain.)
The {\em dimension} of a Noetherian space $Y$ is the maximum length of
a chain of closed irreducible subsets (so the dimensions $\dim \kk$ and
$\dim \Max (\kk)$ of $\kk$ and, respectively, $\Max (\kk)$ are the
same).

\begin{lemma}\label{lemma.nilp}\ 
Let $A$ be a finite type $\kk$ algebra. Then $\Jac(A)$ is
nilpotent. In particular, if $J_1 \subset J_2$ are two--sided ideals
of $A$ such that $V(J_1)=V(J_2)$, then $J_2/J_1$ is nilpotent.
\end{lemma}

\begin{proof}\ 
The result is well known if $A$ is commutative (it is Hilbert's
Nullstellensatz, see \cite{AtiyahMacDonald1,Bourbaki}). We begin with
two remarks.

Let $a \in A$ be an arbitrary element. Assume that $a$ is not
nilpotent, then there exists a maximal ideal $\mfk m \subset \kk[a]$,
$a \not \in \mfk m$. We shall denote by $\chi : \kk[a] \to \CC$ the
unique morphism with kernel $\mfk m$. Thus $\chi(a) \not = 0$.
Then the left $A$-module
\begin{equation*}
	V := A \otimes_{\kk[a]} (\kk[a]/\mfk m)
\end{equation*}
is a finite dimensional complex vector space and hence the vector $\xi
= 1 \otimes 1 \in V$ is contained in an irreducible left $A$-submodule
$W \subset V$. By construction, we have that
$$
	a(1 \otimes 1) = a \otimes 1 = 1 \otimes a = 
	\chi(a) 1 \otimes 1,
$$
so the element $a$ will not act on $W$ by a nilpotent endomorphism.
Therefore $a$ is not in the Jacobson radical of $A$. This shows that
$\Jac(A)$ consists entirely of nilpotent elements.
 
Let $\KK$ be the ring of fractions of $\kk$ (the ring of fractions
associated to the multiplicative subset of elements that are not
zero-divisors). Denote by $Tor(A) \subset A$ the set of torsion
elements of $A$, that is, the set of elements $a \in A$ such that
there exists $f \in \kk$, $f$ not a divisor of zero in $\kk$, but $f a =
0$. By definition, $Tor(A)$ is the kernel of the map $A \to A
\otimes_\kk \KK$, and hence $Tor(A)$ is a two-sided ideal of $A$ and a
$\kk$-submodule. Since $A \otimes_{\kk} \KK$ is a finite dimensional
algebra over a field, $\Jac(A \otimes_{\kk} \KK)$ is nilpotent
(Engel's theorem, see \cite{Fulton}, Theorem 9.9).  Since $\Jac(A) \otimes_\kk \KK$
consists of nilpotent elements, we also have that $\Jac(A) \otimes_\kk
\KK \subset \Jac(A \otimes_\kk \KK)$, and hence $\Jac(A)^N \subset
Tor(A)$, for some large $N$.

We shall prove the lemma by induction on the dimension of the support
of $A$ as a $\kk$-module. By replacing $\kk$ with a quotient, we can
assume that the support of $A$ is the whole of $\Max(\kk)$.  When the
dimension of $\Max(\kk)$ is zero, $\kk$ is actually finitely
dimensional, and the result again follows from Engel's theorem.

Since $\kk$ is noetherian and $A$ is finitely generated, the
submodule $Tor(A)$ is also finitely generated.  Hence we can find $f
\in \kk$, $f$ not a zero divisor in $\kk$, such that $f Tor(A) =
0$. Consequently, by a standard result (see \cite{AtiyahMacDonald1},
Corollary 11.9), the support of $Tor(A)$ has smaller dimension than
the support of $A$. But
\begin{equation*}
	\Jac(A) \cap Tor(A) \subset \Jac(Tor(A)),
\end{equation*}
and hence $\Jac(A)^N \subset \Jac(Tor(A))$.  This shows that the
induction hypothesis can be applied to $Tor(A)$.
  
The last part follows because $J_2/J_1$ is contained in the Jacobson
radical of the finite type algebra $A/J_1$.
\end{proof}

We shall need a few simple facts from the theory of PI-algebras.
Recall \cite{Rowen1} that a {\em polynomial identity algebra} (or,
simply, a {\em PI--algebra}) is an algebra $A$ for which there exists
a non-zero polynomial $P(X_1,X_2, \ldots, X_m)$ in the non-commuting
variables $X_1,X_2, \ldots, X_m$, without constant term, such that
\begin{equation*}
	P(a_1, a_2, \ldots, a_m)=0,
\end{equation*}
for any $a_1, a_2, \ldots, a_m \in A$. The polynomial $P$ is called an
{\em identity} of $A$.  A finite type $\kk$--algebra $A$ satisfies an
alternating multilinear identity of degree $q+1$, where $q$ is the
number of elements in a system of generators of $A$ as a
$\kk$--module, so finite type algebras are PI-algebras.

The {\em PI--class} of a finite type algebra $A$ is, by definition,
the smallest integer $n$ such that any identity of $M_n(\CC)$ is also
an identity of $A$. The PI--class of $A$ has representation theoretic
significance because the largest integer $n$ such that $\prn$ is not
empty coincides with the PI--class of $A/\Jac(A)$.

We shall need a few basic results involving Azumaya algebras, whose definition 
we recall below.

\begin{definition} \
An {\em Azumaya algebra} {\em [of rank $r$]} over the commutative
algebra $Z$ is a unital $Z$--algebra $A$ which is a finitely generated
projective module {\em [of rank $r$]} over $Z$ such that
\begin{equation*}
	A \otimes_{Z} A^{\opp} \cong \End{Z}A.
\end{equation*}
\end{definition}

It follows that $Z$ coincides with the center of $A$.  In general, the
rank of $A$ can be different on different connected components of the
maximal ideal spectrum of its center. Assume that $A$ has rank
$r$. The dimension $[A_{\mfk p}:Z_{\mfk p}]$ of the localization
$A_{\mfk p}=(Z\smallsetminus {\mathfrak p})^{-1}A$ is equal to $r$, for
all maximal ideals $\mfk p \subset A$, and $r=n^2$, where $n$ is the
PI--class of $A$. Moreover, the central character defines a
homeomorphism $\Theta:\PR A \to \PR Z$, see \cite{Pierce1}.  Azumaya
algebras will play an important r\^ole in what follows, because of the
following basic result.

\begin{theorem}[Artin-Procesi]\label{Artin.Procesi} \ 
Suppose that $\mfk A$ is a unital PI--algebra of PI--class $n$ and
that $\PR {\mfk A} =\prim {n}{\mfk A}$, then $\mfk A$ is an Azumaya
algebra of rank $n^2$.
\end{theorem}

See \cite[1.8.48]{Rowen1} or the original papers \cite{Artin1,Procesi}
for a proof.

We shall now prove two lemmas on finite type algebras. We shall assume
in the two lemmas below (and in a few other lemmas later on) that
$\kk$ is an {\em integral domain}, that is, that $\kk$ has no
zero-divisors. We shall do so even when the results hold in greater
generality, because this level of generality is enough for us and the
proofs are shorter. Moreover, the proof of our final results can be
reduced, from the general case of a finitely generated ring $\kk$, to
the case when $\kk$ is an integral domain.

If $\pp \subset \kk$ is a prime ideal and $M$ is a
$\kk$-module, then as usual $M_\pp$ denotes the localization of $M$ at
$\pp$.

\begin{lemma}\label{lemma.Sone}\ Let $A$ be a semi-primitive, finite 
type $\kk$-algebra, and assume that $\kk$ is an integral domain. Then
the set
\begin{equation*}
	S_1 : = \{\pp \in \Max(\kk), \ A_\pp \text{ \rm\ is unital}
	\,\}
\end{equation*}
is open and dense in $\Max(\kk)$.
\end{lemma}

We agree that the algebra $\{0\}$ {\em does have} a unit.

\begin{proof}\  
To prove that $S_1$ is dense, it is enough to show that it is a
non-empty open set, because $\kk$ is an integral domain.

Denote by $\KK$ the ring of fractions of $\kk$. Then $A \otimes_\kk
\KK$ is a semi-primitive finite dimensional $\KK$-algebra, and hence
it has unit.  Let $a \in A$ and $f \in \kk^*$ be such that $f^{-1}a$
is the unit of $A \otimes_\kk \KK$. Let $b_1,\ldots,b_l$ be generators
of $A$ as a $\kk$-module.  Then there exists $g_1,\ldots,g_l \in
\kk^*$ such that $g_j(ab_j - fb_j) = 0$ in $A$, for any $j = 1,
\ldots, l$. If we then let $S$ be the principal open set associated
to $h = fg_1\ldots g_l \not = 0$ (that is, the set of maximal ideals
of $\kk$ not containing $h$), we get $S \subset S_1$. This shows that
$S_1$ is not empty.

To prove that $S_1$ is open, we proceed in an analogous manner.  Let
$\pp$ be a maximal ideal such that $A_\pp$ is unital. Let $f^{-1}a$
(where $a \in A$ and $f \in \kk^*$, $f \not \in \pp$) be that
unit. Define $g_1,\ldots,g_l \in \kk^*$ and $h = fg_1\ldots g_l \not =
0$ as above. Then the principal open set associated to $h$ contains
$\pp$ and is in turn contained in $S_1$. This completes the proof.
\end{proof}

The following lemma is a well known consequence of the
Artin-Procesi theorem. We include a proof for the convenience of the
reader. Let us recall first that a polynomial $g$ in non-commuting
variables is called a {\em central polynomial for $M_n(\CC)$} if, by definition, 
$g(T_1,\ldots,T_r)$ is a multiple of the identity matrix for any choice
of the matrices $T_1, \ldots, T_r \in M_n(\CC)$. Central polynomials are
known to exist for any $M_n(\CC)$, see \cite{Formanek}.

\begin{lemma}\label{lemma.Stwo}\ 
Let $A$ be a unital, semi-primitive, finite type $\kk$-algebra, and
assume that $\kk$ is an integral domain. Then the set
\begin{equation*}
	S_2 : = \{\pp \in \Max(\kk), \ A_\pp \text{ \rm\ is an Azumaya
	algebra}\,\}
\end{equation*}
is open and dense in $\Max(\kk)$. 
\end{lemma}

\begin{proof}\  
Let $Z=Z(A)$ be the center of $A$. By standard commutative algebra,
the set of maximal ideals $\pp \in \Max(\kk)$ such that $A_\pp$ is a
projective $Z_\pp$-module is an open subset of $\Max(\kk)$.  Since
$A$ is a finitely generated $Z$-module, $[\End{Z}A \big ]_\pp =
\End{Z_\pp}{A_\pp}$. Then, the set of maximal ideals $\pp \in \Max(\kk)$ such
that
$$
	A_\pp \otimes_{Z_\pp} A_\pp = (A \otimes_Z A)_\pp \cong \big
	[ \End{Z}A \big ]_\pp = \End{Z_\pp}{A_\pp}
$$ 
is an open subset of $\Max(\kk)$. This shows that $S_2$ is open. To
prove that $S_2$ is dense, it is enough then to prove that it is
not empty.

Let $n$ be the $PI$-type of $A$. Then $\prim{n}A \not = 0$, because
$A$ is semi-primitive. Let $\mfk I \subset A$ be the intersection of
the kernels of all irreducible representations of $A$ of dimension $<
n$. Then $Z$, the center of $A$, intersects $\mfk I$ non-trivially.
Indeed, this is a deep fact that is proved using central polynomials as
follows:\ if $g_n$ is a central polynomial for $M_n(\CC)$, then
$g_n(\mfk I) \subset \mfk I \cap Z$ and $g_n(\mfk I) \not = \{0\}$.)
Let $I := Z \cap \mfk I$, then an irreducible complex representation of
$A$ is in $\prim{n}A$ if, and only if, it does not vanish on $I$. Let
$\mfk q$ be a maximal ideal of $\kk$ {\em not} containing $\kk \cap
\mfk I$. Then $A_{\mfk q}$ is a $PI$-algebra of $PI$-type $n$ and all
irreducible representations of $A_{\mfk q}$ are of dimension $n$. The
Artin-Procesi theorem then asserts that $A_{\mfk q}$ is an Azumaya
algebra. This implies that $S_2$ is not empty and hence completes
the proof.
\end{proof}

\section{Spectrum preserving morphisms}

We shall now study a class of morphisms implicitly appearing in
Lusztig's work on the representation of Iwahori-Hecke algebras, see
\cite{LC3}. More precisely, we shall study the class of ``spectrum
preserving'' morphisms of finite type algebras. Our main reason for
studying this class of morphisms is that such a morphism induces
an isomorphism in periodic cyclic homology (Theorem~\ref{theorem.HPiso}).

Denote by $\Theta_L : \Prim(L) \to \Max(\kk)$ and by $\Theta_J :
\Prim(J) \to \Max(\kk)$ the central character maps of the finite type
$\kk$-algebras $L$ and $J$, Equation \eqref{inf.character2}. If $\phi
: L \to J$ is a $\kk$-linear morphism, we define
\begin{equation}\label{eq.sp.rel}
	{\mathcal R}_\phi := \{ (\mfk P',\mfk P) \subset \Prim(J)
	\times \Prim(L),\, \phi^{-1}(\mfk P') \subset \mfk P \}.
\end{equation}

\begin{lemma}\label{lemma.eq.cc}\
If $(\mfk P', \mfk P) \in {\mathcal R}_\phi$, then $\Theta_L(\mfk P) =
\Theta_J(\mfk P')$.
\end{lemma}

\begin{proof}\ 
If $\mfk P$ and $\mfk P'$ are as in the statement of the lemma, then
an irreducible representation with kernel $\mfk P$ identifies with a
subspace of an irreducible representation with kernel $\mfk P'$. This
shows that any $a \in \kk$ will act on both spaces via the same scalar
$\chi(a)$. Consequently, $\Theta_L(\mfk P) := \ker(\chi)
=:\Theta_J(\mfk P')$.
\end{proof}

We shall need the following lemma.

\begin{lemma}\label{lemma.comp}\ 
Suppose $\phi : L \to J$ and $\psi : J \to A$ are morphisms of finite
type $\kk$-algebras.  Then ${\mathcal R}_{\psi \circ \phi} = {\mathcal
R}_{\psi} \circ {\mathcal R}_{\phi}$.
\end{lemma}

\begin{proof}\ 
The inclusion ${\mathcal R}_{\psi \circ \phi} \supset {\mathcal
R}_{\psi} \circ {\mathcal R}_{\phi}$ follows directly from the
definition.

To prove the other inclusion, suppose $(\pi,V_\pi)$ is an irreducible
representation of $L$ (on the vector space $V_\pi$) with kernel $\mfk
P$.  Let $\mfk P'' \subset A$ be the kernel of an irreducible
representation $(\pi'', V_{\pi''})$ of $A$. If $(\psi \circ
\phi)^{-1}(\mfk P'') \subset \mfk P$, then $V_\pi$ identifies with a
subrepresentation of $\pi'' \circ \psi \circ \phi$. Restrict $\pi''$
to $J$, and let $W$ be an irreducible $J$-module containing $V_\pi$.
Then the annihilator $\mfk P'$ of $W$ satisfies $\phi^{-1}(\mfk P')
\subset \mfk P$ and $\psi^{-1}(\mfk P'') \subset \mfk P'$.  This
proves that ${\mathcal R}_{\psi \circ \phi} \subset {\mathcal
R}_{\psi} \circ {\mathcal R}_{\phi}$.
\end{proof}

We now introduce the class of morphisms we are interested in.

\begin{definition}\label{Def.mult.one}\
Let $\phi : L \to J$ be a $\kk$-linear morphism of unital, finite type
$\kk$-algebras. We say that $\phi$ is a {\em spectrum preserving
morphism} if, and only if, the set $\mathcal R_\phi$ defined in
Equation \eqref{eq.sp.rel} is the graph of a bijective function
\begin{equation*}
	\phi^* : \Prim(J) \to \Prim(L).
\end{equation*}
\end{definition}

More concretely, we see that $\phi : L \to J$ is spectrum preserving if,
and only if, the following two conditions are satisfied:

\begin{enumerate}
\item\ {\em For any primitive ideal $\mfk P$ of $J$, the ideal
$\phi^{-1}(\mfk P)$ is contained in a unique primitive ideal of $L$,
namely $\phi^*(\mfk P)$}, and
\item\ {\em The resulting map $\phi^*: \Prim(J) \to \Prim(L)$ is a
bijection}.
\end{enumerate}

It is instructive to look at the case when both $L$ and $J$ are finite
dimensional, semi-simple. Then any morphism $\phi: L \to J$ gives rise
to a morphism 
\begin{equation}
	\phi_* : K_0(L) \to K_0(J),
\end{equation}
and the conjugacy class of $\phi$ is determined by $\phi_*$. By
identifying the $K_0$-groups with the free abelian subgroups generated
by the simple factors, we see that $\phi_*$ is determined by a matrix
with integer entries, called the {\em Bratteli diagram} of $\phi$
\cite{Bratteli} (see also \cite{Dadarlat}).  A morphism $\phi : L \to
J$ of finite dimensional semi-simple algebras will then be spectrum
preserving if, and only if, it decomposes up to conjugacy as a direct
sum of (not necessarily unital) injective morphisms $M_r(\CC) \to
M_q(\CC)$. This condition is equivalent to saying that the Bratteli
diagram of $\phi$ is a diagonal matrix with non-zero diagonal entries
(after reordering the simple factors of $J$, if necessary). The
definition of the Bratteli diagram extends also to other classes of
rings $L$ and $J$ in the obvious way. We shall be interested below
(Lemma \ref{lemma.widehat}) in the case when $L$ and $J$ are direct
sums of matrix algebras over some completion of $\kk$.

{}From Lemma~\ref{lemma.comp} we get the following proposition, which
implies, among other things, that the composition of two spectrum
preserving morphisms is again spectrum preserving.

\begin{proposition}\label{prop.comp}\ 
Suppose that $\phi : L \to J$ and $\psi : J \to A$ are morphisms of
finite type $\kk$-algebras. If two of the morphisms $\phi$, $\psi$, or
$\psi \circ \phi$ are spectrum preserving, then the third morphism is
spectrum preserving as well and $(\psi \circ \phi)^* = \phi^* \circ
\psi^*$.
\end{proposition}

\begin{proof}\ 
The relation ${\mathcal R}_{\psi \circ \phi} = {\mathcal R}_{\psi}
\circ {\mathcal R}_{\phi}$ of the above lemma tells us right away that
if two of the relations ${\mathcal R}_\phi$, ${\mathcal R}_\psi$, or
${\mathcal R}_{\psi \circ \phi}$ are the graphs of bijective
functions, then the third one is also the graph of a bijective
function.  Moreover, the relation $(\psi \circ \phi)^* = \phi^* \circ
\psi^*$ also follows.
\end{proof}

Assuming that the morphism $\phi : L \to J$ is spectrum preserving, we
obtain the following description of $\phi^*$.  If
$\mfk P$ is a primitive ideal of $J$, denote by $p : L \to L /
\phi^{-1}(\mfk P)$ the natural projection. Then
$$
	\phi^*(\mfk P) = p^{-1}\big [\Jac(L / \phi^{-1}(\mfk P)) \big ].
$$
An explicit description of $(\phi^*)^{-1}$ will be obtained later on.

It is easy to check using the definition that $\phi: L \to J$ is a
spectrum preserving morphism if, and only if, $\phi^+ : L^+ \to J^+$
is a spectrum preserving morphism, and hence we can take the latter to
be the definition of a spectrum preserving morphism in the non-unital
case.  It also follows directly from the definition that $L \to J$ is
spectrum preserving if, and only if, $L \to J/\Jac(J)$ is spectrum
preserving.

\begin{lemma}\label{lemma.inj.J}\
The kernel of a spectrum preserving morphism $\phi: L \to J$ is
contained in the Jacobson radical of $L$. Consequently, a morphism
$\phi : L \to J$ is spectrum preserving if, and only if, $\ker(\phi)
\subset \Jac(L)$ and the induced map $L/\ker(\phi) \to J$ is spectrum
preserving.
\end{lemma}

\begin{proof}\ Let $\mfk Q \subset L$ be an arbitrary primitive ideal. 
We need to show that $\ker(\phi) \subset \mfk Q$. By the surjectivity
of the map $\phi^* : \Prim(J) \to \Prim(L)$, there exists a primitive
ideal $\mfk P \subset J$ such that $\phi^{-1}(\mfk P) \subset \mfk Q$.
Then $\ker(\phi) : = \phi^{-1}(0) \subset \mfk Q$.
\end{proof}

{}From the above discussion and the above lemma, we get the following.

\begin{corollary}\label{cor.inj.J}\ 
A morphism $\phi : L \to J$ of finite type algebras is spectrum
preserving if, and only if, $\phi^{-1}(\Jac(J)) \subset \Jac(L)$ and
the induced morphism 
\begin{equation*}
	L/\phi^{-1}(\Jac(J)) \to J /\Jac(J)
\end{equation*}
is spectrum preserving.
\end{corollary}

Spectrum preserving morphisms behave well when we restrict them to
ideals and when we take quotients.

\begin{proposition}\label{prop.sp.pres}\ 
Let $\phi : L \to J$ be a $\kk$-linear morphism of (not necessarily
unital) finite type algebras, and let $J' \subset J$ be a two-sided
ideal. If $\phi$ is spectrum preserving, then the induced morphisms 
$\phi^{-1}(J') \to J'$ and $L / \phi^{-1}(J') \to J / J'$ are spectrum 
preserving.
\end{proposition}

\begin{proof}\ 
By considering algebras with adjoined unit, we can assume that both
$L$ and $J$ are unital. The Jacobson spectrum $\Prim(J)$ decomposes as
a disjoint union of the set of primitive ideals containing $J'$ and
the set of primitive ideals {\em not} containing $J'$. These two sets
identify with the Jacobson spectra of $J/J'$ and $J'$, respectively.

By Lemma~\ref{lemma.inj.J}, we can assume that $\phi$ is an inclusion
(so $L$ identifies with a subalgebra of $J$). We can restrict the
morphism $\phi^*$ to each of $\Prim(J/J')$ and $\Prim(J')$, the two
sets considered above.

It is clear from definitions that $\phi^*$ maps a primitive ideal of
$J$ containing $J'$ to a primitive ideal of $L$ containing $L \cap
J'$. This means that $\phi^*$ maps $\Prim(J/J')$ to $\Prim(L/(L \cap
J'))$. It is then enough to show that the induced map
$$
	\phi^*: \Prim(J/J') \to \Prim(L/(L \cap J'))
$$ 
is a bijection.

Lemma~\ref{lemma.eq.cc} gives that $\phi^*$ induces a bijection of the
sets $\Theta_J^{-1}(\pp)$ and $\Theta_L^{-1}(\pp)$.  We also know that
$\phi^*$ maps
\begin{equation*}
	S := \Theta_J^{-1}(\pp) \cap \Prim(J/J') \;\; \text{ to } \; \;
	S' := \Theta_L^{-1}(\pp) \cap \Prim(L/(L \cap J')).
\end{equation*}
To prove that $\phi^*$ maps $\Prim(J/J')$ to $\Prim(L/(L \cap J'))$
bijectively, it is enough to prove that $\phi^* : S \to S'$ is
bijective.

Now we have that $S = \Prim(J/\pp J)$ and $S' = \Prim(L/\pp L)$.
Moreover, it is clear from Lemma~\ref{lemma.eq.cc} that the induced
map $\phi_0 : L / \pp L \to J / \pp J$ is spectrum preserving and that
$\phi_0^* = \phi^*$ as maps from $S$ to $S'$.  Since both $J/\pp J$
and $L /\pp L$ are finite dimensional, this reduces our discussion to
the finite dimensional case.

So assume now that $L$ and $J$ are finite dimensional, that $J$ is
semi-primitive, and that $L \subset J$.  Then $J$ is semi-simple.  
Since the exact sequence
\begin{equation}
	0 \to \Jac(L) \to L \to L / \Jac(L) \to 0
\end{equation}
splits, we can also assume that $L$ is semi-primitive.  If we write
$J= \oplus J_\alpha$ as a direct sum of simple factors, the assumption
that $L \to J$ is spectrum preserving means that we can write $L =
\oplus \, (L \cap J_\alpha)$, from which the proposition follows.
\end{proof}

It is interesting to mention the following consequence of the (proof of the)
above proposition.

\begin{corollary}\label{cor.V}\ 
If $\phi : L \to J$ is spectrum preserving and $J' \subset J$ is an
ideal, then $\phi^*(V(J')) = V(\phi^{-1}(J'))$.
\end{corollary}

Similarly, we have the following property.

\begin{lemma}\ Assume that $\phi : L \to J$ is spectrum preserving and 
that all irreducible representations of $J$ have dimension $ \le n $,
then $(\phi^*)^{-1}(V(L')) = V(J \phi(L')^n J)$, for any two-sided
ideal $L' \subset L$.
\end{lemma}

\begin{proof}\ 
Let $\mfk P \subset J $ be an arbitrary primitive ideal. Since
$L/\phi^{-1}(\mfk P)$ identifies with a subalgebra of the algebra of
$m \times m$ matrices for some $m \le n$, the Jacobson radical of $
L/\phi^{-1}(\mfk P) $ is nilpotent of order at most $m$, by Engel's
theorem. Hence, if $L_1:= \Jac(L/\phi^{-1}(\mfk P) )$, then $L_1^n =
0$. This shows that $\phi^*(\mfk P)$, which is the inverse image of
$L_1$ in $L$, satisfies $\phi^*(\mfk P)^n \subset \phi^{-1}(\mfk P)$.

We have $\mfk P \in (\phi^*)^{-1}(V(L'))$ if, and only if, $L' \subset
\phi^*(\mfk P)$. The above discussion shows that $L' \subset
\phi^*(\mfk P) \Leftrightarrow {L'}^n \subset \phi^{-1}(\mfk
P)$. Since
\begin{equation*}
	{L'}^n \subset \phi^{-1}(\mfk P) \Leftrightarrow \phi(L')^n
	\subset \mfk P \Leftrightarrow J \phi(L')^n J \subset \mfk P,
\end{equation*}
the result follows.
\end{proof}

By applying the above lemma to a primitive ideal $\mfk P'$ of $L$, we
obtain the following description of $(\phi^*)^{-1}$. Denote by $p : J
\to J/J \phi(\mfk P') J$ the natural projection, and let $n$ be as in
the above lemma. Then
$$
	(\phi^*)^{-1}(\mfk P') = p^{-1} ( \Jac( J/J \phi(\mfk P')^n J
	)).
$$
Together, the above results give the following theorem.

\begin{theorem}\ 
Let $L$ and $J$ be finite type algebras and $\phi : L \to J$ be a
spectrum preserving morphism. Then the induced map $\phi^* : \Prim(J)
\to \Prim(L)$ is a homeomorphism.
\end{theorem}

The following result, which is, in a way, a converse of
Corollary~\ref{cor.V}, will be needed later on.

\begin{corollary}\label{cor.pres.V}\ 
Suppose $\phi : L \to J$ is a spectrum preserving morphism and let
$L_1 \subset L$ and $J_1 \subset J$ be two-sided ideals such that
$\phi(L_1) \subset J_1$.  Then the induced map $L_1 \to J_1$ is
spectrum preserving if, and only if, $\phi^*(V(J_1)) = V(L_1)$. Similarly, 
the induced map $L/L_1 \to J/J_1$ is spectrum preserving if, and only
if, $\phi^*(V(J_1)) = V(L_1)$.
\end{corollary}

\begin{proof}\ We know that $L_2:=\phi^{-1}(J_1) \to J_1$ is spectrum 
preserving (by Proposition~\ref{prop.sp.pres}) and that $V(L_2) =
\phi^*(V(J_1))$. Since $L_1 \subset L_2$ is a two-sided ideal of
$L_2$, the inclusion $L_1 \to L_2$ is spectrum preserving if, and only
if, $V(L_1) = V(L_2)$.  The fact that the map $L_1 \to J_1$ is
spectrum preserving then follows from Proposition~\ref{prop.comp}.

Let us prove that $L/L_1 \to J/J_1$ is spectrum preserving.  Indeed,
Corollary \ref{cor.V} and the assumption that $\phi^*(V(J_1)) =
V(L_1)$ show that $L/L_1 \to L/\phi^{-1}(J_1)$ is spectrum
preserving. The corollary then follows from Proposition \ref{prop.comp}
and Proposition \ref{prop.sp.pres}.
\end{proof}

\section{Periodic cyclic homology \label{Sec.cyclic}}

In this section, we recall a few basic results on cyclic and
Hochschild homology necessary for the proof, in the next section, of
the fact that a spectrum preserving morphism induces an isomorphism in
periodic cyclic homology. The account of Hochschild and cyclic
homology in \cite{KNS} is more complete, and we refer the reader
without prior exposure to cyclic homology to that paper for more
homological results relevant to our setting. See also
\cite{ConnesDG,ConnesNCG,Karoubi,LodayQuillen}, and \cite{Tsygan}.

Let $A$ be a complex unital algebra.  Define the operators
$s,t,B,b',b$, acting on $A^{\otimes n + 1} : = A \otimes A \otimes
\ldots \otimes A$ ($n+1$ times), as follows.
\begin{equation*}
%\label{equation.b}
\begin{gathered}
	s(\Tt) = 1\otimes \Tt, \\
	t(\Tt) = (-1)^n a_n\otimes a_0\otimes\ldots\otimes a_{n-1}, \\
 	B(\Tt) = (1-t)s\sum_{k=0}^{n} t^k(\Tt),\\
	b'(\Tt) = \sum_{i=0}^{n-1} (-1)^ia_0\otimes\ldots\otimes 
	a_i a_{i+1}\otimes\ldots\otimes a_n,\;\quad \text{ and }\\
	b(\Tt)=b'(\Tt) + (-1)^n a_na_0\otimes\ldots\otimes a_{n-1}. 
\end{gathered}
\end{equation*}
Here we have used the notation of \cite{ConnesNCG}. Then the {\em
Hochschild homology groups of $A$}, denoted $\Hd_q(A)$, can be
computed as the homology groups of the complex $(A^{\otimes
n+1},b)$. If $A$ is a $\kk$-algebra, then $\Hd_q(A)$ is a
$\kk$-module.  The {\em cyclic homology} groups of $A,$ denoted
$\Hc_n(A),$ are the homology groups of the cyclic complex $({\mathcal
C}(A),b+B)$, where
\begin{equation*}
	{\mathcal C}_n(A) = \bigoplus_{k\geq 0}A\otimes(A/\CC
	1)^{\otimes n-2k}.
\end{equation*}
The complex $\mathcal C_n$ has an endomorphism $S : \mathcal C_n \to
\mathcal C_{n-2}$. The inverse limit of this complex with respect to
the morphism $S$ is called the
{\em periodic cyclic homology complex}. Its homology is called the {\em
periodic cyclic homology} of $A$ and is denoted by $\Hp_*(A)$.

Let $J$ be a two--sided ideal of $A$. As usual, 
$\hat A := {\lim_\leftarrow} A/J^k$. By replacing in the definitions
above $A^{\otimes n+1}$ with
\begin{equation*}
	A^{\Wh\otimes n+1} := \displaystyle{\lim_\leftarrow}\,
	(A/J^k)^{\otimes n+1}, \quad k \to \infty,
\end{equation*}
we obtain the definitions of the groups $\tHd_*(\Wh A)$ and
$\tHp_*(\Wh A)$. We shall use these constructions for $A$ a
$\kk$-algebra, where $\kk$ is a Noetherian ring, and $J = IA$, for
some ideal $I$ of $\kk$.  Denote, as usual by
\begin{equation*}
	\Wh{M} := \displaystyle{\lim_{\leftarrow}}\, M/I^n M,
	\quad n \to \infty,
\end{equation*}
the $I$-adic completion of a $\kk$--module $M$.  If $M$ is a finitely
generated $\kk$-module, then $\Wh{M}\cong M \otimes_{\kk}\ckk$, see
\cite{AtiyahMacDonald1,Bourbaki}.

The following result \cite[Theorem 3]{KNS} describes the effect of
$I$--adic completion on Hochschild homology.

\begin{theorem}\label{theorem.completion}\
Assume that $\kk_0$ is a noetherian ring and that $A_0$ is a unital
$\kk_0$--algebra, which is a finitely generated $\kk_0$-module. Let $I
\subset \kk_0$ be an ideal, and let $\Wh{A}_0 :=
\displaystyle{\lim_{\leftarrow}}\, A_0/ I^k A_0$ be the completion of
$A_0$ in the $I$-adic topology. Then the natural map $\Hd_*(A_0) \to
\tHd_*(\Wh{A}_0)$ and the $\kk_0$--module structure on
$\Hd_*(A_0)$ define an isomorphism
\begin{equation*}
	\Hd_*(A_0) \otimes_{\kk_0} \Wh{\kk}_0 \cong
	\tHd_*(\Wh{A}_0)
\end{equation*} 
of $\Wh \kk_0$--modules.
\end{theorem}

In \cite{KNS}, the above theorem was proved when $\kk_0$ a finitely
generated algebra, however, the proof works without any change in the
case of a noetherian ring $\kk_0$. For the convenience of the reader,
we include a proof of this theorem in the case when $\kk_0 = \kk_\pp$
is the localization of a finitely generated commutative ring at a
maximal ideal $\pp$ and $I = \pp\kk_\pp$, which is the only case we
shall need in this paper.

\begin{proof}\ 
The theorem will be proved by applying Theorem 3 of \cite{KNS} and
Proposition \ref{Prop.localization}, which is a result from
\cite{Brylinski1}.

We know that the theorem is true when $\kk_0$ is a finitely generated
ring (Theorem 3 from \cite{KNS}). Assume now that $\kk$ is a finitely
generated commutative ring and that $\pp$ is a maximal ideal of
$\kk$. Let $\kk_0 = \kk_{\pp}$, $A_0 = A_{\pp}$, and $I = \pp
A_{\pp}$. We shall reduce the proof of our theorem in this case to a
straightforward application of Theorem 3 of \cite{KNS}. According
to that theorem
\begin{equation*}
	\tHd_*(\Wh{A}) \cong \Hd_*(A) \otimes_{\kk} \Wh{\kk}.
\end{equation*}
Therefore,
\begin{equation*}
	\tHd_*(\Wh{A}_0) = \tHd_*(\Wh{A}) \cong \Hd_*(A) \otimes_{\kk}
	\Wh{\kk} = \Hd_*(A) \otimes_{\kk} \Wh{\kk_0}.
\end{equation*}
On the other hand, by Proposition~\ref{Prop.localization},
\begin{equation*}
	\Hd_*(A_0) \cong \Hd_*(A) \otimes_\kk \kk_0.
\end{equation*}
Hence,
\begin{equation*}
	\Hd_*(A_0) \otimes_{\kk_0} \Wh{\kk}_0 \cong (\Hd_*(A)
	\otimes_\kk \kk_0) \otimes_{\kk_0} \Wh{\kk}_0 \cong \Hd_*(A)
	\otimes_\kk \Wh{\kk}_0.
\end{equation*}
Thus $\Hd_*(A_0) \otimes_{\kk_0} \Wh{\kk}_0 \cong \tHd_*(\Wh{A}_0)$.
\end{proof}

Recall now the proposition from \cite{Brylinski1} used above.

\begin{proposition}\label{Prop.localization} \
Let $S$ be a multiplicative subset of the center $Z$ of the algebra
$A$. Then $\Hd_*(S^{-1}A) \cong S^{-1}\Hd_*(A)$.
\end{proposition}

The following result is well known (a proof of it can be found in
\cite{KNS}, for example). A proof can also be obtained from Proposition
\ref{Prop.localization}.

\begin{proposition}\label{etale.triv} \
Let $Z$ be a commutative ring and $A$ be an Azumaya algebra over $Z$.
Then the inclusion $j:Z \to A$ defines an isomorphism $j_*:\Hd_*(Z)
\to \Hd_*(A)$.
\end{proposition}

We now recall two basic results on cyclic homology. We begin with the
following result of Goodwillie \cite{Goodwillie1}.

\begin{theorem}\label{Goodwillie} \
If $I \subset A$ is a nilpotent ideal, then the quotient morphism $A
\to A/I$ induces an isomorphism $\Hp_*(A) \to \Hp_*(A/I)$.
\end{theorem}

A deep result of Cuntz and Quillen, which we shall often use below, is
the Excision property in periodic cyclic homology
\cite{CuntzQuillen0}.

\begin{theorem}[Excision]\label{Theorem.Excision}
Any two-sided ideal $J$ of an algebra $A$ over a characteristic zero
field gives rise to a natural, six-term periodic exact sequence

\begin{equation}
\xymatrix{
\Hp_0(J)        		\ar[r] &
\Hp_0(A)			\ar[r] & 
\Hp_0(A/J) 			\ar[d]^{\pa} \\
\Hp_1(A/J) 			\ar[u]^{\pa} &  
\Hp_1(A)        		\ar[l] &
\Hp_1(J).        		\ar[l] }
\end{equation} 
\end{theorem}

An immediate consequence is that $\Hp_*(I) = 0$ if $I$ is a nilpotent 
algebra.

\section{A comparison theorem}

We shall prove now our main result on cyclic homology,
Theorem~\ref{theorem.HPiso}. Theorem~\ref{theorem.HPiso} formulates a
general principle, which, in rough terms, states that a morphism of
finite type algebras that preserves the spectrum induces an
isomorphism of periodic cyclic homology groups.  Thus the periodic
cyclic homology of a finite type algebra is a ``spectral invariant.''

We begin with a sequence of lemmas. Recall that we denote by $M_\pp =
S^{-1}M$ the localization of a $\kk$-module $M$ with respect to a
maximal ideal $\pp \subset \kk$. Denote by $\Wh{\kk}$ the
completion of $\kk$ with respect to (the powers of) $\pp$.

\begin{lemma}\label{lemma.widehat}\ 
Assume that $\Wh{L}$ and $\Wh{J}$ are Azumaya algebras and that their
centers are isomorphic to finite direct sums of copies of
$\Wh{\kk}$. Let $\Wh{\phi} : \Wh{L} \to \Wh{J}$ be a $\Wh{\kk}$-linear
algebra morphism. If the induced morphism $\phi_1 : \Wh{L}/\pp \Wh{L}
\to \Wh{J}/\pp \Wh{J}$ is spectrum preserving, then $\Wh{\phi}_* :
\tHd_q(\Wh{L}) \to \tHd_q(\Wh{J})$ is an isomorphism for all $q$.
\end{lemma}

\begin{proof}\ Assume that the center of $\Wh{L}$ is isomorphic to $\ckk^m$. 
Then 
\begin{equation*}
	\tHd_*(\Wh{L}) \cong \tHd_*(\ckk)^m \cong \ckk^m,
\end{equation*}
each direct summand corresponding to a
direct summand of the center of $\Wh{L}$. The morphism
$\Wh{\phi}_* : \tHd_q(\Wh{L}) \to \tHd_q(\Wh{J})$ identifies then
with a matrix, which is immediately seen to be isomorphic to the
Bratteli diagram of the morphism $\phi_1$.
\end{proof}

We shall now use Lemma \ref{lemma.widehat} to prove Lemma 
\ref{lemma.Sthree}.

\begin{lemma}\label{lemma.Sthree}\  
Let $\phi : L \to J$ be a unital, injective, spectrum preserving
morphism of semi-primitive, finite type $\kk$-algebras. Assume that
$\kk$ is an integral domain. Then the set
\begin{equation*}
	S_3 : = \{\pp \in \Max(\kk), \ \phi_*: \Hd_*(L_\pp) \to
	\Hd_*(J_\pp) \, \text{ is an isomorphism }\}
\end{equation*}
is open and dense in $\Max(\kk)$.
\end{lemma}

\begin{proof}\ For any morphism $f : M \to M'$ of finitely generated
$\kk$-modules, it is known that the set
\begin{equation*}
	\{\pp \in \Max(\kk),\ f_{\pp}: M_\pp \cong M_\pp'\,\}
\end{equation*}
is open in $\Max(\kk)$, the maximal ideal spectrum of $\kk$.  Both 
$\Hd_*(L)$ and $\Hd_*(J)$ are finitely generated $\kk$-modules
(see \cite{KNS} for a proof). This observation and
Proposition~\ref{etale.triv} imply that $S_3$ is open.

To prove that $S_3$ is dense, it is enough to prove that it is not
empty, because $\Max(\kk)$ is irreducible.  Since $L$ and $J$ are
unital and semi-primitive, using Lemma~\ref{lemma.eq.cc} we get that
we can also assume $\kk$ to map injectively to each of $L$ and $J$.
By Lemma~\ref{lemma.Stwo}, we know that we can find an open dense set
$S_2 \subset \Max(\kk)$ such that both $L_\pp$ and $J_\pp$ are Azumaya
algebras.

Let $Z(L)$ be the center of $L$ and $\pi_L : \Max(Z(L)) \to \Max(\kk)$
be the central character map defined by the inclusion $\kk \to Z(L)$,
which is defined because $Z(L)$ is a finite type $\kk$-algebra. Define
$\pi_J : \Max(Z(J)) \to \Max(\kk)$ similarly.

Let $\widetilde \kk$ be the field of fractions of the integral domain
$\kk$, and let $n_L$ be the dimension of $Z(L) \otimes_{\kk}
\widetilde \kk$ as a $\widetilde \kk$-vector space.  Define $n_J$
similarly. We can choose a maximal ideal $\pp \in S_2$ such that the
following three conditions are satisfied:\ $Z(L)/\pp Z(L)$ has
dimension $n_L$, $Z(J)/\pp Z(J)$ has dimension $n_J$, and all points
of $\pi_L^{-1}(\pp)$ and $\pi_J^{-1}(\pp)$ are regular. (Recall that a
maximal ideal $\mfk m \subset \kk$ is called {\em regular} if $\dim
\,(\mfk m/\mfk m^2) = \dim \kk$, see \cite{AtiyahMacDonald1}, Theorem
11.22.  Geometrically, this means that the associated point in the
spectrum of $\kk$ is non-singular.)  We shall prove then that $\pp \in
S_3$.

Let $\Wh{\kk}$ be the completion of $\kk$ with respect to $\pp$.
Since $\ckk$ is faithfully flat over $\kk_\pp$
(\cite{AtiyahMacDonald1}, Theorem 10.17) it is enough to prove
that the map $\phi$ induces an isomorphism
\begin{equation*}
	\Hd_*(L_\pp)\otimes_\kk \ckk \cong \Hd_*(J_\pp) \otimes_{\kk}
	\ckk.
\end{equation*}
Let $\Wh{L}$ and $\Wh{J}$ be the completions of $L$ and $J$ with
respect to $\pp$. Using that $L_\pp \otimes_{\kk}\ckk \cong \Wh L$
and $J_\pp \otimes_{\kk}\ckk \cong \Wh J$, we obtain from
Theorem~\ref{theorem.completion} that it is enough to prove that
\begin{equation*}
	\tHd_*(\Wh{L}) \cong \tHd_*(\Wh{J}).
\end{equation*}
By assumptions, $\Wh{L}$ is an Azumaya algebra over $\Wh{Z(L)} := Z(L)
\otimes_{\kk} \ckk$. Since $\pp$ is regular, the ring $\ckk$ is
isomorphic to the ring of formal power series.  Similarly, because the
preimage of $\pp$ in $\Max(Z(L))$ consists of regular points, the ring
$\Wh{Z(L)}$ is a direct sum of rings of power series.  More precisely,
we have that $\Wh{Z(L)} \cong \ckk^m$, where $m$ is the degree of the
field of fractions of $Z(L)$ over the field of fractions of $\kk$.
The result then follows from Lemma \ref{lemma.widehat} and Lemma
\ref{lemma.eq.cc}.
\end{proof}

We shall also need the following lemma from \cite{KNS}

\begin{lemma}\label{lemma.third} \
Let $L $ be a $\kk$--algebra and $\pp \subset \kk$ be a maximal
ideal such that $L_{\pp}$ is unital. If $B = \kk +L$ is the
$\kk$--algebra with adjoined unit, then the inclusions $L \to B$ and
$\kk \to B$ induce isomorphisms $\Hd_*(B)_{\pp} \cong
\Hd_*(\kk_{\pp}) \oplus \Hd_*(L_{\pp})$.
\end{lemma}

In particular, $\Hd_*(L_\pp) = \ker [ \Hd_*(B)_{\pp} \to
\Hd_*(\kk)_{\pp}]$. The following proposition puts together the above
results.

\begin{proposition}\label{proposition.Sfour}\  
Let $\phi : L \to J$ be an injective spectrum preserving morphism
of semi-primitive, finite type $\kk$-algebras, with $\kk$ an integral
domain. Then the set
\begin{equation*}
	S_4 : = \{\pp \in \Max(\kk), \  L_\pp \text{ \rm\ and } J_\pp
	\text{ \rm\ are unital and } \Hd_*(L_\pp) \cong \Hd_*(J_\pp) \,\}
\end{equation*}
is open and dense in $\Max(\kk)$.  
\end{proposition}

\begin{proof}\ We know that the set 
\begin{equation*}
	S_3 := \{\pp \in \Max(\kk), \ \phi_* : \Hd_*(L^+_\pp) \cong
	\Hd_*(J^+_\pp) \,\}
\end{equation*}
is open and dense in $\Max(\kk)$ by Lemma~\ref{lemma.Sthree}.  By
Lemma~\ref{lemma.Sone}, we know that the set
\begin{equation*}
	S_1 : = \{\pp \in \Max(\kk), \   L_\pp \text{ \rm\ and } J_{\pp} 
	\text{ \rm\ are unital}\, \}
\end{equation*}
is also open and dense in $\Max(\kk)$. Then Lemma~\ref{lemma.third}
gives that $S_4 = S_1 \cap S_3$, from which the result follows.
\end{proof}

We shall need also one of the main results from \cite{KNS}, Theorem 8
(or rather Corollary 5 of that paper) which gives a criterion for two
algebras to have the same periodic cyclic homology.

\begin{theorem}[Kazhdan-Nistor-Schneider]\label{theorem.Mult}\ 
Let $A$ be a unital finite type $\kk$--algebra and $L \subset J
\subset IA$ be inclusions of finite type $\kk$--algebras, where $I
\subset \kk$ is an ideal. Suppose that, for all maximal ideals
${\mathfrak p}$, $I \not \subset \mfk p$, the localizations
$L_{\mathfrak p}$ and $J_{\mathfrak p}$ are unital and
$\Hd_q(L_{\mathfrak p}) \to \Hd_q(J_{\mathfrak p})$ is an
isomorphism. Then the inclusion $L \to J$ gives an isomorphism
$$\Hp_q(L) \cong \Hp_q(J).$$
\end{theorem}

We are ready now to prove the main result on cyclic homology of this paper,
a result that reflects the {\em spectral invariance
of periodic cyclic homology} in the class of finite type
algebras.

\begin{theorem}\label{theorem.HPiso}\ Let $\phi : L \to J$ be a spectrum 
preserving morphism of finite type $\kk$-algebras. Then the induced map
$\phi_* : \Hp_*(L) \to \Hp_*(J)$ is an isomorphism.
\end{theorem}

\begin{proof}\ 
By Lemma~\ref{lemma.eq.cc}, the $\kk$-modules $L$ and $J$ have the
same support
\begin{equation*}
	Y:=\supp L = \supp J, 
\end{equation*}
a closed subset of $\Max(\kk)$ in the Zariski topology.  We can
replace $\kk$ with a quotient and assume that $Y = \Max(\kk)$.

We shall prove the theorem by induction on $n = \dim \Max(\kk)$.  If
$\dim \Max(\kk) = 0$, then $\kk$ is finite dimensional and hence $J$
and $L$ are also finite dimensional algebras. We can then use the
argument at the end of the proof of Proposition \ref{prop.sp.pres} to
replace $L$ and $J$ by their semi-simple quotients, in which case the
result is straightforward.

We now proceed with the inductive step. Assume that
\begin{quote}
{\em $(A)$\ if $\phi_0 : L_0 \to J_0$ is a spectrum preserving
morphism of finite type $\kk_0$-algebras and $\dim \kk_0 < n$,
then $(\phi_0)_* : \Hp_*(L_0) \to \Hp_*(J_0)$ is an isomorphism.}
\end{quote}
We want to prove then that $\phi_* : \Hp_*(L) \to \Hp_*(J)$ is an
isomorphism for all spectral preserving morphisms between 
finite type $\kk$-algebras with $\dim \kk = n$.

We shall use the induction assumption and excision in periodic cyclic
homology to successively reduce the proof of the isomorphism $\phi_* :
\Hp_*(L) \to \Hp_*(J)$ to the proof of an isomorphism of the same
kind, but with different algebras $L$ and $J$, which can then be
handled directly. We now begin our sequence of (six) reductions.

\vspace*{1mm}{\em Step 1.}\ First, let $\operatorname{nil}(\kk)$ be
the nilradical of $\kk$.  Since $\operatorname{nil}(\kk)L$ and
$\operatorname{nil}(\kk)J$ are nilpotent, they are contained in the
Jacobson radicals of $L$ and, respectively, $J$ and have vanishing
periodic cyclic homology groups. By Corollary \ref{cor.pres.V}, the
induced morphism $\phi_1 : L/\operatorname{nil}(\kk)L \to
J/\operatorname{nil}(\kk)J$ is also spectrum preserving. By excision
and Lemma \ref{lemma.nilp}, $\phi_*$ is an isomorphism if, and only
if, $(\phi_1)_*$ is an isomorphism.  We can then replace $L$ by
$L/\operatorname{nil}(\kk)L$, $J$ by $J/\operatorname{nil}(\kk)J$, and
$\kk$ by $\kk/\operatorname{nil}(\kk)$.  This shows that we can also
assume that $\kk$ is reduced.

\vspace*{1mm}{\em Step 2.}\ Then, using Lemma~\ref{lemma.nilp}, the
fact that a nilpotent ideal has vanishing periodic cyclic homology,
and Lemma~\ref{lemma.inj.J}, we see that, after replacing $L$ by
$L/\ker(\phi)$, if necessary, we can also assume that the morphism
$\phi$ is injective.

\vspace*{1mm}{\em Step 3.}\ Using induction on the number of
irreducible components of $\Max(\kk)$, we can further reduce to the
case when $\kk$ is an integral domain, as follows. Let $\pp_1, \pp_2,
\ldots, \pp_r$ be the minimal prime ideals of $\kk$. Then $\pp_1 \cap
\pp_2 \cap \ldots \cap \pp_r = 0,$ because we have assumed $\kk$ to be
reduced. Let
\begin{equation*}
	J_j := (\pp_1 \cap \ldots \cap \pp_j) J
\end{equation*}
and $L_j := \phi^{-1}(J_j)$. Then the induced maps $L_j/L_{j + 1} \to
J_j/J_{j + 1}$ are injective and spectrum preserving, by
Proposition~\ref{prop.sp.pres}.  Moreover, the $\kk$-module structure
on $J_j/J_{j + 1}$ descends to a $\kk/\pp_{j+1}$-module
structure. This shows that, after replacing $\kk$ by $\kk/\pp_{j+1}$
and $L \to J$ by $L_j/L_{j + 1} \to J_j/J_{j + 1}$, we can also assume
that $\kk$ is an integral domain.

\vspace*{1mm}{\em Step 4.}\ Using Lemma~\ref{lemma.nilp}, the fact
that a nilpotent ideal has vanishing periodic cyclic homology, and
Corollary~\ref{cor.inj.J}, we can also assume that $J$ is
semi-primitive, eventually after replacing $J$ by $J/\Jac(J)$ and $L$
by $L/\phi^{-1}(\Jac(J))$.

\vspace*{1mm}{\em Step 5.}\ For any $\kk$-module $M$, we shall denote
by $Tor(M)$ the set of torsion elements of $M$, that is, the set of
elements whose annihilator contains at least one non-zero element. We
already know that $Tor(L)$ and $Tor(J)$ are two-sided ideals in $L$
and, respectively, $J$ (see the proof of
Lemma~\ref{lemma.nilp}). Since $\phi : L \to J$ is injective,
$Tor(L)=\phi^{-1}(Tor(J))$. Proposition \ref{prop.sp.pres} then implies
that both $Tor(L) \to Tor(J)$ and $L/Tor(L) \to J/Tor(J)$ are
injective and spectrum preserving. Due to the fact that the support of $Tor(L)$ (=
the support of $Tor(J)$, by Lemma \ref{lemma.eq.cc}) has smaller Krull
dimension, using the excision property of periodic cyclic homology,
we see from the inductive hypothesis that it is enough to prove that
$L/Tor(L) \to J/Tor(J)$ induces an isomorphism in periodic cyclic
homology. In other words, we can also assume that both $L$ and $J$
are torsion free.

We need to check however that this reduction does not affect the
previous reductions. More precisely, we need to check that the
quotient $J/Tor(J)$ is still semi-primitive. Indeed, $\Jac(J/Tor(J))$
coincides with the nil-radical of $J/Tor(J)$, so it is enough to prove
that if $I \subset J$ is a two-sided ideal such that $I^n \subset
Tor(J)$, then $I \subset Tor(J)$. If $I^n \subset Tor(J)$, because
$\kk$ is noetherian, there exists $f \not = 0$ such that $fI^n =
0$. But then $(fI)^n = 0$, so $fI = 0$, because $J$ is semi-primitive.

\vspace*{1mm}{\em Step 6.}\ Our next reduction is to show that we can
assume in addition that $L$ is semi-primitive.  Let $\KK = S^{-1}\kk$ be the
field of fractions of $\kk$ (so $S = \kk^*$ because we have reduced to
the case when $\kk$ is an integral domain). The algebra $B := L
\otimes_{\kk} \KK$ is a finite dimensional algebra over the field
$\KK$ and it satisfies
\begin{equation*}
	\Jac(B) \cong \Jac(L) \otimes_{\kk} \KK \quad \text{and} \quad
	B/\Jac(B) \cong (L/\Jac(L)) \otimes_{\kk} \KK.
\end{equation*} 
Moreover, $B$ splits, in the sense that $B$ contains a subalgebra
which maps isomorphically to $B/Jac(B)$. Since $L$ is torsion free,
$L/\Jac(L)$ is also torsion free, and hence $B$ will also contain an
algebra $L_1$ which projects isomorphically onto $L/\Jac(L)$ via the
natural map $L_1 \to B \to B/\Jac(B)$.

Since $L_1$ is a finitely generated $\kk$-module, there exists then
$f \in \kk$, $f \not = 0$, such that $fL_1 \subset L$ (we used here
the fact that $L$ is torsion free, so it identifies with a subalgebra
of $B$).

Let $J_1 := f J$. Then $\phi : L \to J$ maps $L_2:=f^2L_1 \subset f L$
to $J_1$. We have that $V(L_2) = \Theta_L^{-1}(V(f^2\kk))$ and $V(J_1)
= \Theta_J^{-1}(V(f\kk))$.  Consequently,
\begin{equation*}
	\phi^*(V(J_1)) = V(L_2),
\end{equation*}
and hence the induced map $L_2 \to J_1$ is spectrum preserving, by
Corollary~\ref{cor.pres.V}. By construction, both $L_2$ and $J_1$ are
semi-primitive.

Moreover, the $\kk$-module structures on $L/L_2$ and $J/J_1$ descend
to $\kk/(f^2)$-module structures, and the maximal ideal spectrum of
$\kk/(f^2)$ has dimension $(\dim \kk) - 1$. Using the excision
property of periodic cyclic homology and the inductive hypothesis, we
see that it is enough to prove that $L_2 \to J_1$ induces an
isomorphism in periodic cyclic homology.  We obtain finally that in
the induction step we can also assume both $L$ and $J$ to be
semi-primitive.  This completes our sequence of
reductions. \vspace*{1mm}

We stress that each reduction is in addition to the previous
reductions.  In this way, we have obtained that -- in order to
establish the inductive step -- it is enough to prove the following:

\begin{quote}{\em $(R)$:\ if $\phi : L \to J$ is a spectrum 
preserving morphism of finite type $\kk$-algebras, where $\kk$ is an
integral domain, $\phi$ is injective, and $L$ and $J$ are
semi-primitive, then $\phi_* : \Hp_*(L) \to \Hp_*(J)$ is an
isomorphism.}
\end{quote}
(We no longer need $L$ or $J$ to be torsion free as $\kk$-modules, that
assumption was necessary only to reduce to the case where we can
assume $L$ to be semi-primitive.)

We now prove $(R)$, assuming $(A)$ (the inductive assumption).  The
advantage of $(R)$ is that we can now use
Proposition~\ref{proposition.Sfour} to conclude that the set
\begin{equation*}
	S_4 : = \{\pp \in \Max(\kk), \ L_\pp \text{ \rm\ and } J_\pp
	\text{ \rm\ are unital and } \phi_* : \Hd_*(L_\pp) \cong
	\Hd_*(J_\pp) \,\}
\end{equation*}
is open and dense in $\Max(\kk)$. Let $I \subset \kk$ be the ideal
corresponding to this open set. Define $L_1 : = L \cap IJ^+$ and $J_1
:= J \cap IJ^+ = IJ$.  Then $L_1 \subset J_1 \subset IJ^+$ satisfy the
assumptions of Theorem ~\ref{theorem.Mult} for the ideal $I$ (we use
here also that $I_\pp \cong \kk_\pp$, if $I \not \subset \pp$), and
hence the inclusion $L_1 \to J_1$ induces an isomorphism $\Hp_*(L_1)
\to \Hp_*(J_1)$.

The induced morphism $\phi_1:L/L_1 \to J/J_1$ is spectrum preserving,
by Proposition~\ref{prop.sp.pres}. By construction, the $\kk$-module
structure on the quotients $L/L_1$ and $J/J_1$ descends to a
$\kk/I$-module structure, and the maximal ideal spectrum of $\kk/I$
has smaller dimension than the dimension of $\Max(\kk)$.  By the
induction hypothesis, $\phi_1$ induces an isomorphism in periodic
cyclic homology. The excision property then gives that $L \to J$
induces also an isomorphism in periodic cyclic homology. The induction
step is now established, so the theorem is proved.
\end{proof}

\section{Periodic cyclic homology of Iwahori-Hecke algebras\label{Sec.Hecke}}

Since the converse of Proposition~\ref{prop.sp.pres} is not true, the
assumption that a morphism be spectrum preserving is usually too
restrictive in applications. However, the following weaker version of
the spectrum preserving property is usually enough.

\begin{definition}\label{def.w.sp.p}\ 
A morphism $\phi : L \to J$ of finite type algebras is called {\em
weakly spectrum preserving} if, and only if, there exist increasing
filtrations
\begin{equation*}
\begin{gathered}
	(0) = L_0 \subset L_1 \subset \ldots \subset L_n = L \quad
	\text{ and } \qquad (0) = J_0 \subset J_1 \subset \ldots
	\subset J_n = J
\end{gathered}
\end{equation*} 
of two-sided ideals such that $\phi(L_k) \subset J_k$ and the induced
morphisms $L_k/L_{k-1} \to J_k / J_{k-1}$ are spectrum preserving.
\end{definition}

Using the excision property of periodic cyclic homology, we see that
Theorem~\ref{theorem.HPiso} immediately implies the following only
slightly stronger, but much more useful, form of that theorem:

\begin{theorem}\label{theorem.HPiso2}\ 
Let $\phi : L \to J$ be a weakly spectrum preserving morphism of
finite type $\kk$-algebras. Then the induced map $\phi_* : \Hp_*(L)
\to \Hp_*(J)$ is an isomorphism.
\end{theorem}

It is interesting to note that if we denote by $R(A)$ the Grothendieck
group of finite dimensional modules over $A$ and if $\phi: L \to J$ is
a weakly spectrum preserving morphism, then we obtain an induced map
$\phi^* : R(J)\otimes \QQ \to R(L) \otimes \QQ$, which is an
isomorphism. (Compare to \cite{LC3} and note that in that case there
is no need to tensor with $\QQ$.) Also, we obtain bijections
$\Prim(J_k/J_{k-1}) \to \Prim(L_k/L_{k-1})$, which in turn induce a
bijection $\Prim(J) \to \Prim(L)$. This bijection is not continuous,
in general, and it may also depend on the choice of the filtrations
used to check the weakly spectrum preserving property.

We now explain how the above theorem and the results of
Kazhdan-Lusztig \cite{KL} and Lusztig \cite{LC1,LC2,LC3,LC4} and
\cite{Asterisque} lead to a determination of the periodic cyclic
homology of Iwahori-Hecke algebras $\Hecke_q$, for $q \in \CC^*$ not a
proper root of unity. (See also \cite{Springer} for $q=1$.)  We
begin by recalling the classical definition of Hecke algebras (we
shall work here with affine Hecke algebras and extended affine Hecke
algebras).

Fix a root system $\Phi \subset E$ in a Euclidean space and let $W$ be
the affine Weyl group generated by $\Phi$, that is, the group of
affine isometries of $E$ generated by the reflections with respect to
all hyperplanes $H_{\alpha,k} := \{ \lambda \in E, \, (\lambda,\alpha)
= k\,\}$, where $k \in \ZZ$ and $\alpha \in \Phi$. (See
\cite{Humphreys}).  Let $W_0$ be the finite Weyl group associated to
$\Phi$ (it is the group generated by the reflections relative to the
planes $H_{\alpha,0}$. It is well known that $W$ is a Coxeter group
(see \cite{Humphreys}, Theorem 4.6). We fix a system of generators $S
\subset W$ such that $(W,S)$ is a Coxeter system.  The resulting
length function on $W$ will be denoted by $l$, as in
\cite{Asterisque}.

Then the {\em Iwahori-Hecke algebra $\Hecke_q$ associated to $W$} and $q \in \CC^*$ 
is the algebra generated by $T_{x}$, $x \in W$, with relations
\begin{equation}\label{eq.def.Hecke}
\begin{gathered}
	T_xT_y = T_{xy}, \;\; \text{if}\; l(xy) = l(x) + l(y), \quad
	\text{and}\\ (T_s - q)(T_s + 1) = 0 , \quad \text{if} \; s \in
	S.
\end{gathered}
\end{equation} 
As a vector space over $\CC$, the algebra $\Hecke_q$ has a basis consisting
of $T_x$, $x \in W$.

Let $\Omega$ be a group acting by automorphisms on $(W,S)$. That is,
$\Omega$ acts by automorphisms on $W$ and maps $S$ to itself.  Let
$\hW : = W \rtimes \Omega$. We extend $l$ to $\hW$ by
$l(w\omega)=l(w)$, if $w \in W$ and $\omega \in \Omega$.  Then the
{\em Iwahori-Hecke algebra $\Hecke_q$ associated to $\hW$} is the
algebra generated by $T_{x}$, $x \in \hW$, and with exactly the same
relations as those of the Equation \ref{eq.def.Hecke}.

We now proceed to recall the definition of Lusztig's ``asymptotic
Iwahori-Hecke algebra'' $J$. We need first to recall the definition of
the generic Hecke algebra $\gHecke_\param$ and of the Kazhdan-Lusztig
polynomials. Recall that $\gHecke_\param$ is the (complex) algebra
with generators $\param$, $r$, and $T_x$, $x \in W$, satisfying the
relations
\begin{equation}\label{eq.def.gHecke}
\begin{gathered}
	T_x T_y = T_{xy}, \;\; \text{if}\; l(xy) = l(x) + l(y), \qquad
	rT_x = T_xr,\; r^2 = \param,\,  \qquad \text{and}\\ (T_s -
	\param)(T_s + 1) = 0 , \quad \text{if} \; s \in S.
\end{gathered}
\end{equation} 
The map $r \to r^{-1}$ and $T_w \to T_{w^{-1}}^{-1}$, $w \in W$,
defines an involution $\kappa : \gHecke_\param \to
\gHecke_\param$. Then \cite{KL0}, for any $w \in W$, there exists a
unique element $C_w \in \gHecke_\param$ such that $\kappa(C_w) = C_w$
and
\begin{equation}\label{eq.def.KLp}
	C_w = r^{-l(w)} \sum_{y \le w} P_{y,w}T_y, \quad y \in W,
\end{equation}
where $P_{w,w} = 1$ and $P_{y,w}$ is a polynomial in $r$ of degree
$\le l(w) - l(y) -1$, if $y < w$. The polynomials $P_{y,w}$ are called
the {\em Kazhdan-Lusztig} polynomials.

Let $h_{w,u,v} \in \CC[r,r^{-1}]$ be defined by
\begin{equation}
	C_wC_u= \sum_{v \in W} h_{w,u,v}C_v.
\end{equation}
For any $v \in W$, let $a(v)$ be the smallest non-negative integer $i$
with the property that $r^i h_{w,u,v}$ has no strictly negative powers
of $r$.  Then denote by $\gamma_{w,u,v}$ the constant term of $(-r)^i
h_{w,u,v}$.

The {\em asymptotic Hecke algebra} $J$ is then defined as the algebra
generated by $t_w$, $w \in W$, and relations
\begin{equation*}
	t_x t_y = \sum_{z \in W} \gamma_{x,y,z}t_z.
\end{equation*}

Let ${\mathcal D} := \{ w \in W, \, \deg P_{e,w} = l(w) - a(w)\,\}.$
Define the map
\begin{equation}\label{eq.Lusztig}
	\phi : \gHecke_\param \to J \otimes \CC[r,r^{-1}], \quad
	\phi(C_w) = \sum_{a(d) = a(u)} h_{w,d,u} t_u,
\end{equation}
where in the above sum $w \in W$ and $d \in {\mathcal D}$ are
arbitrary subject to the condition that $a(d) = a(u)$. Then $\phi$ is
a morphism of $\CC[r,r^{-1}]$ algebras, see \cite{LC2} (see also
\cite{NanhuaXi}, page 23). Denote by $\phi_q : \Hecke_q \to J$ any of
the morphism obtained by specializing at $\param = q$ and by choosing
a square root of $q$ (which induces a morphism $\CC[r,r^{-1}] \to
\CC$). The resulting morphism depends only on $q$ and not on the
particular choice of a square root of $q$.

Fix now $q \in \CC^*$. We know that if we denote by $\kk$ the center
of $\Hecke_q$, then $\Hecke_q$ and $J$ are finite type $\kk$-algebras,
see \cite{LC3}, Proposition~1.6 (the $\kk$-module structure on $J$ is
obtained via $\phi_q$).

\begin{theorem}[Lusztig]\ If $q  \in \CC^*$ is not a proper root of unity,
then the morphism $\phi_q : L = \Hecke_q \to J$ is weakly spectrum preserving.
\end{theorem}

\begin{proof}\  
Let $J^{(i)} \subset J$ be the ideal generated by $t_w$, $a(w)=i$.
Then $J^{(i)} \subset J$ is a two-sided ideal and $J \cong \oplus_i
J^{(i)}$. Define $J_k = \oplus_{i \ge k} J^{(i)}$ and $L_k \subset
L:=\Hecke_q$ by $L_k = \phi_q(J_k)$. Then the results of
\cite{Asterisque}, Theorem 8.1, prove the theorem.
\end{proof}

{\em Remark.}\ The filtration $J_k$ used in the above proof should be
compared with the filtrations from \cite{LC3}, Theorem 3.4 and
Corollary 3.6, or \cite{NanhuaXi}, 5.9.  The result in
\cite{Asterisque} does not involve taking the ``dual'' $E \to E^*$ of
a $J$-module, so we get our chosen inequality sign, the other two
references provide the opposite inequality sign, that is, we have to
define then $J_k = \oplus_{i \le k} J^{(i)}$.\vspace*{.8mm}

A consequence of the above computation is the following result on the
periodic cyclic homology of the algebras $J$ and $\Hecke_q$, for $q
\in \CC^*$ not a proper root of unity.

\begin{theorem}\ 
Assume that $q \in \CC^*$ is not a proper root of unity. Then
Lusztig's morphism $\phi_q : \Hecke_q \to J$, see Equation
\eqref{eq.Lusztig}, induces an isomorphism
\begin{equation*}
	(\phi_q)_*: \Hp_*(\Hecke_q) \to \Hp_*(J).
\end{equation*}
\end{theorem}

For $q=1$ we can find $\Hp_*(\Hecke_q)$ explicitly, so this leads to a
determination of the groups $\Hp_*(\Hecke_q)$, provided that $q \in
\CC^*$ is not a proper root of unity. Indeed, let $\hW$ be the
extended affine Weyl group used in the definition of the algebras
$\Hecke_q$ and $J$. Then $\Hecke_1 \cong \CC[\hW]$.

Assume now that $\hW$ acts properly and isometrically on a Euclidean
space $E$ (this is the case for example if $\Omega$ is
commutative). The space $E$ contains the root system $\Phi$, but is
not necessarily generated by $\Phi$.  It follows then that any
subgroup $G \subset \hW$ will have finite cohomological dimension over
a characteristic zero field. Indeed, this is because
\begin{equation}
	\cohom_*(G;\CC) \cong \cohom_*(E/G;\CC).
\end{equation}

We continue to assume that $\hW$ acts properly and isometrically on a
Euclidean space $E$. Denote by $o(g)$ the linear part of an isometry
$g$, defined by writing the group of isometries as the direct product
of the group of linear isometries (=isometries fixing the origin) with
the group of translations.  If $g \in G$, we shall denote by $\hW_g :=
\{ g'\in \hW, \, gg' = g'g\, \}$, that is the {\em centralizer of $g$
in $\hW$}.  If $g \in \hW$ is an element of infinite order, let $E_0$
be the subspace fixed by $o(g)$. If $g' \in G$ commutes with $g$, then
$g'$ maps $E_0$ to itself. By definition, $g$ acts by translations on
$E_0$. Let $E'$ be the one-dimensional subspace generated by the
direction of this translation and let $(g) \cong \ZZ$ be the subgroup
generated by $g$. Then $E_0/E'$ is acted upon {\em properly} by
$\hW/(g)$. This shows that all groups of the form $\hW_g/(g)$ have
finite cohomological dimension.

Let $\langle G \rangle'$ be the set of finite order conjugacy classes
of a group $G$, let $\langle G \rangle''$ be the set of infinite order
conjugacy classes of $G$ and $\langle G \rangle := \langle G \rangle'
\cup \langle G \rangle''$.  If $g \in G$, we shall denote by $G_g :=
\{ g'\in G, \, gg' = g'g\, \}$, that is the {\em centralizer of $g$ in
$G$}.  Recall now a result of Burghelea \cite{Burghelea} (see also
\cite{Karoubi} and \cite{Ni1}), which states that the {\em cyclic
homology groups} of $\CC[G]$ are given by
\begin{equation}\label{eq.Burghelea}
	\Hc_p(\CC[G]) \cong \bigoplus_{\langle g \rangle \in \langle G
	\rangle'} \left ( \oplus_{k \ge 0} \cohom_{p - 2k}(G_g;\CC)
	\right ) \bigoplus \bigoplus_{\langle g \rangle \in \langle G
	\rangle''} \cohom_p(G_g/(g);\CC) .
\end{equation}
Moreover, the natural morphism $S : \Hc_{p+2}(\CC[G]) \to
\Hc_{p}(\CC[G])$ preserves the above direct sum decomposition
according to conjugacy classes and is the natural projection on each
component corresponding to a finite order conjugacy class.

This gives the following description of the groups $\Hp_p(\Hecke_q)$
and $\Hp_*(\CC[\hW])$.

\begin{theorem}\ 
Assume that $q \in \CC^*$ is not a proper root of unity and let $\hW$
be an extended affine Weyl group acting properly on a Euclidean space
$E$. Then the map $S : \Hc_{p+2}(\CC[\hW]) \to \Hc_{p}(\CC[\hW])$ is
surjective for $p \ge \dim E$ and
\begin{multline*}
	\Hp_p(\Hecke_q) \cong \Hp_p(\CC[\hW]) \cong \bigoplus_{\langle
	g \rangle \in \langle \hW \rangle'} \left ( \oplus_{k \in \ZZ}
	\cohom_{p - 2k}(\hW_g;\CC) \right ) \\ \cong
	\bigoplus_{\langle g \rangle \in \langle \hW \rangle'} \left (
	\oplus_{k \in \ZZ} \cohom_{q - 2k}(E/\hat W_g;\CC) \right ).
\end{multline*}
All these isomorphisms are natural.
\end{theorem}

We note that it is sometimes convenient to use $\cohom_*(E/\hW_g)
\cong \cohom_*(E^g/\hW_g)$, if $g$ is a finite order element of $\hat
W$. We conclude this section with a few remarks.

Let $\Omega$ be a discrete group acting by automorphisms on the finite
type $\kk$--algebras $L$ and $J$, and let $L \rtimes \Omega$ and $J
\rtimes \Omega$ be the corresponding crossed-products. If $\phi : L
\to J$ is $\Omega$-equivariant, spectrum preserving morphism, then the
induced map $\phi \rtimes \Omega : L \rtimes \Omega \to J \rtimes
\Omega$ is also spectrum preserving. This is proved using an explicit
description of the primitive ideal spectra of $L \rtimes \Omega $ and
$J \rtimes \Omega $.  If $\phi$ is only weakly spectrum preserving,
but the defining filtrations of Definition~\ref{def.w.sp.p} are
$\Omega$-invariant, then $\phi \rtimes \Omega : L \rtimes \Omega \to J
\rtimes \Omega $ will be again weakly spectrum preserving. The example
we have in mind is when $L$ is the Iwahori-Hecke algebra associated to
an affine Weyl group $W$ and $J$ is the corresponding asymptotic Hecke
algebra.  Then $L \rtimes \Omega$ and $J \rtimes \Omega$ are the
Hecke, respectively, asymptotic Hecke algebras of $\hW$, for a
suitable action of $\Omega$.

%%%%%%%%  BIBLIOGRAPHY  %%%%%%%%

\providecommand{\bysame}{\leavevmode\hbox to3em{\hrulefill}\thinspace}


\begin{thebibliography}{10}

\bibitem{Artin1} M.~Artin, \emph{On {A}zumaya algebras and finite
dimensional representations of rings}, J. Algebra \textbf{11} (1969),
532--563.

\bibitem{AtiyahMacDonald1}
M.F. Atiyah and MacDonald, \emph{Introduction to commutative algebra},
  Addison-Wesley, Reading, Mass.-London, 1969.

\bibitem{BC} P.~Baum and A.~Connes,
\emph{Chern character for discrete groups}, A F{\^e}te of Topology,
163--232, Academic Press, New York, 1988.

\bibitem{BBN} P.~Baum, J.~Brodzki, and V.~Nistor, \emph{Cyclic
homology of crossed products}, work in progress.

\bibitem{BCH} P.~Baum, A.~Connes, and N.~Higson, \emph{Classifying
spaces for proper actions and $K$-theory of group $C^*$-algebras},
$C^*$-algebras: 1943--1993 (San Antonio, TX, 1993), 240--291,
Contemp. Math. \textbf{167}, AMS, Providence, RI, 1994.

\bibitem{BHP2} P.~Baum, N.~Higson, and R.~Plymen, \emph{Equivariant
homology for $\operatorname{SL}_2$ of a $p$-adic field}, in Index
theory and operator algebras (Boulder, CO, 1991), 1--18,
Contemp. Math. \textbf{148}, Amer. Math. Soc., Providence, R.I., 1993.

\bibitem{BHP} P.~Baum, N.~Higson, and R.~Plymen, \emph{Proof of the
Baum-Connes conjecture for $p$-adic $GL(n)$}, C. R. Acad. Sci. Paris,
S\'erie I 325 (1997), 171-176.

\bibitem{BHP-hc} P.~Baum, N.~Higson, and R.~Plymen, Representation
theory of $p$-adic groups: a view from operator algebras, \emph{The
mathematical legacy of Harish-Chandra} (Baltimore, MD, 1998), 111-149,
Proc. Symp. Pure Math. 69, A. M. S. Providence RI, 2000.

\bibitem{Borel} A. Borel, \emph{Admissible representations of a
semi-simple group over a local field with vectors fixed under an
Iwahori subgroup}, Invent. Math. \textbf{35} (1976),~233--259.

\bibitem{Bourbaki} N.~Bourbaki, \emph{Alg\`ebre {C}ommutative},
Hermann, Paris, 1961.

\bibitem{Bratteli} O. Bratteli, \emph{Inductive limits of finite
dimensional $C\sp{*} $-algebras}, Trans. Amer. Math. Soc. 171 (1972),
195--234.

\bibitem{Braun} A. Braun, \emph{The nilpotency of the radical in a
finitely generated P.I. ring}, J. Algebra \textbf{89} (1984),
375--396.

\bibitem{Brylinski1} J.-L. Brylinski, \emph{Central localization in
{H}ochschild homology}, J. Pure Appl. Algebra \textbf{57} (1989),
1--4.

\bibitem{Burghelea} D.~Burghelea, 
\emph{The cyclic cohomology of group rings}, 
Comment. Math. Helvetici \textbf{60} (1985), 354--365. 

\bibitem{BushnellKutzko} C. Bushnell and Ph. Kutzko, The admissible dual of
$GL(N)$ via compact open subgroups. Annals of Mathematics
Studies \textbf{129}, Princeton University Press, 1993.

\bibitem{ConnesDG} A.~Connes, \emph{Non commutative differential
geometry}, Publ. Math. IHES \textbf{62} (1985), 41--144.

\bibitem{ConnesNCG} A.~Connes, \emph{Non-commutative geometry},
Academic press, 1994.

\bibitem{CuntzQuillen0} J.~Cuntz and D.~Quillen, \emph{Excision in
bivariant periodic cyclic cohomology}, Invent. Math. \textbf{127}
(1997), 67--98.

\bibitem{Dadarlat} M. Dadarlat, \emph{Nonnuclear subalgebras of AF
algebras}, Amer. J. Math. \textbf{122} (2000), 581--597.

\bibitem{Formanek} E. Formanek, \emph{Central polynomials for matrix
rings}, J. Algebra \textbf{23} (1972), 129--132.

\bibitem{Fulton} W. Fulton and J. Harris, \emph{Representation
Theory}, Graduate Texts in Mathematics \textbf{129}, Springer Verlag,
1991.

\bibitem{GinzburgKV} V. Ginzburg, M. Kapranov, and E. Vasserot,
\emph{Residue construction of Hecke algebras}, Adv. Math.
\textbf{128} (1997), 1--19.

\bibitem{Goodwillie1} T.~G. Goodwillie, \emph{Cyclic homology,
derivations, and the free loopspace}, Topology \textbf{24} (1985),
187--215.

\bibitem{Harris} J. Harris, \emph{Algebraic geometry. A first course},
Grad. Texts in Math. \textbf{133}, Springer-Verlag, 1995.

\bibitem{Higson-Nistor1} N.~Higson and V.~Nistor, \emph{Cyclic
cohomology of groups acting on buildings}, J. Funct.Anal. \textbf{141}
(1996), 466--495.

\bibitem{Humphreys} J. E. Humphreys, \emph{Reflection groups and Coxeter
groups}, Cambridge Studies in Advanced Mathematics 29, Cambridge
University Press, Cambridge, 1990.  xii+204.

\bibitem{IwahoriMatsumoto} N. Iwahori and H. Matsumoto, \emph{On some
Bruhat decomposition and the structure of Hecke rings of $p$-adic
Chevalley groups}, Publ. Math. I.H.E.S.  \textbf{25} (1965), 5--48.

\bibitem{Karoubi} M.~Karoubi, \emph{Homologie cyclique et {K}-theorie}
Asterisque \textbf{149} (1987), 1--147.

\bibitem{KL0} D. Kazhdan and G. Lusztig , \emph{Representations of
Coxeter groups and Hecke algebras}, Invent. Math. \textbf{53} (1979),
165--184.

\bibitem{KL} D. Kazhdan and G. Lusztig , \emph{Proof of the
Deligne-Langlands conjecture for Hecke algebras},
Invent. Math. \textbf{87} (1987), 153--215.

\bibitem{KNS} D.~Kazhdan, V.~Nistor, and P.~Schneider,
\emph{Hochschild and cyclic homology of finite type algebras}, Selecta
Math. (N.S.) \textbf{4} (1998), 321--359.

\bibitem{Kemer} A. T. Kemer, \emph{Capelli identities and nilpotency
of the radical of a finitely generated PI-algebra}, Dokl. Akad. Nauk
SSSR \textbf{255} (1980), 793--797.

\bibitem{Loday} J.-L.~Loday, \emph{Cyclic Homology}, Springer-Verlag,
Berlin-Heidelberg-New York, 1992.
 
\bibitem{LodayQuillen} J.-L. Loday and D.~Quillen, \emph{Cyclic
homology and the {L}ie homology of matrices},
Comment. Math. Helv. \textbf{59} (1984), 565--591.

\bibitem{LC1} G.~Lusztig, \emph{Cells in affine Weyl groups}, In
``Algebraic Groups and Related Topics,'' pp. 255--287, Advanced
studies in Pure Math., vol.  \textbf{6}, Kinokunia and North Holland,
1985.

\bibitem{LC2} G.~Lusztig, \emph{Cells in affine Weyl groups II},
J. Algebra \textbf{109} (1987), 536--548.

\bibitem{LC3} G.~Lusztig, \emph{Cells in affine Weyl groups III},
J. Fac. Sci. Univ. Tokyo Sect. IA Math.  \textbf{34} (1987), 223--243.

\bibitem{LC4} G.~Lusztig, \emph{Cells in affine Weyl groups IV},
J. Fac. Sci. Univ. Tokyo Sect. IA Math.  \textbf{36} (1989), 297--328.

\bibitem{LusztigQG} 
G.~Lusztig,  Introduction to quantum groups, Progress in Mathematics, 110,
Birkh\"auser Boston, Inc., Boston MA, 1993, xii+341.


\bibitem{Asterisque} 
G.~Lusztig,  \emph{Representations of affine Hecke algebras}, 
Ast\'{e}risque \textbf{171-172} (1989), 73--84.  


\bibitem{Manin} Yu. Manin, \emph{Topics in noncommutative geometry},
M. B. Porter Lectures, Princeton University Press, Princeton, NJ,
1991, viii+164pp.

\bibitem{Ni1}
V.~Nistor, \emph{Group cohomology and the cyclic cohomology of crossed
products}, Invent. Math. \textbf{99} (1990), 411--424.

\bibitem{NistorHOMPDC}
V.~Nistor, \emph{Higher orbital integrals, Shalika germs, 
and the Hochschild homology of Hecke algebras}, xxx.math.RT/0008133,
to appear.

\bibitem{Pierce1}
R.~S. Pierce, \emph{Associative algebras}, Grad.
Texts in Math. \textbf{88}, Springer-Verlag, 1982.

\bibitem{Procesi} C.~Procesi, \emph{On a theorem of {M.} {A}rtin},
J. Algebra \textbf{22} (1972), 309--315.

\bibitem{Procesi2} C.~Procesi, \emph{Rings with Polynomial Identities},
Dekker, New York, 1973.

\bibitem{Quillen} D. Quillen, \emph{On the endomorphism ring of a simple
module over an enveloping algebra}, Proc. Amer. Math. Soc. \textbf{21}
(1969), 171--172.

%\bibitem{Rinehart} G.~Rinehart, \emph{Differential forms for general
%commutative algebras}, Trans. Amer. Math. Soc. \textbf{108} (1963),
%195--222.

\bibitem{Rowen1} L.~H. Rowen, \emph{Polynomial identities in ring
theory}, Academic Press, New York-London-Toronto, 1980.

\bibitem{Schneider1} P.~Schneider, \emph{The cyclic homology of
$p$-adic reductive groups}, J. Reine Angew. Math. \textbf{475} (1996),
39--54.

\bibitem{Springer} T. A. Springer, \emph{Linear Algebraic Groups}, 
Birkha\"{u}ser, Boston-Basel-Stuttgart, 1981.

\bibitem{Tsygan} B.~L. Tsygan, \emph{Homology of matrix {L}ie
algebras over rings and {H}ochschild homology}, Uspekhi
Math. Nauk. \textbf{38} (1983), 217--218.

\bibitem{NanhuaXi}
Nanhua Xi, \emph{Representations of affine Hecke algebras}, 
Lecture Notes in Mathematics, 1587, Springer-Verlag, Berlin, 1994, viii+137pp..
\end{thebibliography}
\end{document}